\documentclass{gtart}

\def\ifplaintex{\expandafter\ifx\csname documentclass\endcsname\relax}

\def\gtp{{\mathsurround=0pt\it $\cal G\mskip-2mu$eometry \&\ 
$\cal T\!\!$opology $\cal P\!$ublications}}  % GT publications

\def\recd{{\small Received:\qua\receiveddate\ifx\reviseddate\relax
\else\qquad Revised:\qua\reviseddate\fi\par}} 

\def\lognumber#1{\def\thelognumber{#1}}
\def\volumenumber#1{\def\thevolumenumber{#1}}
\def\volumeyear#1{\def\thevolumeyear{#1}}
\def\papernumber#1{\def\thepapernumber{#1}}
\def\pagenumbers#1#2{\def\startpage{#1}\def\finishpage{#2}}
\def\published#1{\def\publishdate{#1}}

\def\received#1{\def\receiveddate{#1}}
\def\revised#1{\def\reviseddate{#1}}
\def\accepted#1{\def\accepteddate{#1}}

\long\def\asciiabstract#1{\long\def\theasciiabstract{#1}}

\let\\\par\let\thelognumber\relax\let\thevolumenumber\relax
\let\thepapernumber\relax\let\thevolumeyear\relax\let\startpage\relax
\let\finishpage\relax\let\publishdate\relax\let\receiveddate\relax
\let\reviseddate\relax\let\accepteddate\relax\let\theasciititle\relax
\let\theasciiauthors\relax
\let\theasciiabstract\relax

\let\theasciiemail\relax

\ifplaintex
\font\logobig=cmssbx10 scaled 3836
\font\logomed=cmssbx10 scaled 2557
\else
\font\logobig=cmssbx10 scaled 4200
\font\logomed=cmssbx10 scaled 2800
\fi

\long\def\makeagttitle{   %%% start of definition of \makeagttitle
\count0=\startpage
\agt\hfill      %   Journal title (top left) 
\hbox to 45truept{\vbox to 0pt{\vglue -13truept{\logomed A\kern -.37em{\logobig 
T}\kern -.38em G}\vss}\hss}
\break
{\small Volume \thevolumenumber\ (\thevolumeyear)
\startpage--\finishpage\nl
Published: \publishdate}

\vglue .25truein

{\parskip=0pt\leftskip 0pt plus
1fil\def\\{\par\smallskip}{\Large\bf\thetitle}\par\medskip} \vglue
0.05truein

{\parskip=0pt\leftskip 0pt plus 1fil\def\\{\par}{\sc\theauthors}
\par\medskip}%
 
\vglue 0.03truein 

{\small\leftskip 25truept\rightskip 25truept{\bf Abstract}\stdspace\theabstract

{\bf AMS Classification}\stdspace\theprimaryclass
\ifx\thesecondaryclass\relax\else; \thesecondaryclass\fi\par
{\bf Keywords}\stdspace \thekeywords\par}\vglue 7truept

}   %%%% end of definition of \makeagttitle

\ifplaintex
\font\phead=cmsl9 scaled 950
\font\pnum=cmbx10 scaled 913
\font\pfoot=cmsl9 scaled 950
\headline{\vbox to 0pt{\vskip -4.5mm\line{\small\phead\ifnum
\count0=\startpage ISSN 1472-2739 (on-line) 1472-2747 (printed)
\hfill {\pnum\folio}\else\ifodd\count0\def\\{ }% 
\ifx\theshorttitle\relax\thetitle\else\theshorttitle\fi\hfill{\pnum\folio}
\else\def\\{ and }{\pnum\folio}\hfill\ifx\theshortauthors\relax\theauthors
\else\theshortauthors\fi\fi\fi}\vss}}
\footline{\vbox to 0pt{\vglue 0mm\line{\small\pfoot\ifnum\count0=\startpage
\copyright\ \gtp\hfill\else
\agt, Volume \thevolumenumber\ (\thevolumeyear)\hfill\fi}\vss}}
\else
\font=cmsl9 scaled 1050
\font\lnum=cmbx10 
\font=cmsl9 scaled 1050
\makeatletter
\def\@oddhead{{\small\lhead\ifnum\count0=\startpage ISSN 1472-2739 
(on-line) 1472-2747 (printed)\hfill {\lnum\number\count0}\else\ifodd\count0
\def\\{ }\ifx\theshorttitle\relax \thetitle \else\theshorttitle\fi\hfill
{\lnum\number\count0}\else\def\\{ and }{\lnum\number\count0}
\hfill\ifx\theshortauthors\relax 
\theauthors\else\theshortauthors\fi\fi\fi}}\def\@evenhead{\@oddhead}
\def\@oddfoot{\small\lfoot\ifnum\count0=\startpage\copyright\ \gtp\hfill\else
\agt, Volume \thevolumenumber\ (\thevolumeyear)\hfill\fi}
\def\@evenfoot{\@oddfoot}
\makeatother
\fi
\let\maketitlepage\makeagttitle
\let\makeshorttitle\maketitlepage
\let\maketitle\maketitlepage

\newwrite\gtoutfile
\long\gdef\makeheadfile{  %%% start of definition of \makeheadfile
{\def\\{, }\def\s{ }
\immediate\openout\gtoutfile head.xxx
\immediate\write\gtoutfile{To: math@arxiv.org}
\immediate\write\gtoutfile{Subject: put OR rep NNNNN:ppppp}
\immediate\write\gtoutfile{--text follows this line--}
\immediate\write\gtoutfile{Proxy-for: \ifx\theasciiauthors\relax
\theauthors\else\theasciiauthors\fi\s<\ifx\theasciiemail\relax\theemail\else\theasciiemail\fi>}
\immediate\write\gtoutfile{\noexpand\\}
\immediate\write\gtoutfile{Authors: \ifx\theasciiauthors\relax
\theauthors\else\theasciiauthors\fi}
{\def\\{ }\immediate\write\gtoutfile{Title: \ifx\theasciititle\relax
\thetitle\else\theasciititle\fi}}
\immediate\write\gtoutfile{Subj-class: GT or SG, GR etc}
\immediate\write\gtoutfile{MSC-class: \theprimaryclass\ifx\thesecondaryclass\relax\else, \thesecondaryclass\fi}
\immediate\write\gtoutfile{Journal-ref: Algebr. Geom. Topol. \thevolumenumber\s
(\thevolumeyear) \startpage-\finishpage}
\immediate\write\gtoutfile{Comments: Published by Algebraic and
Geometric Topology at}
\immediate\write\gtoutfile{\s\s\s  http://www.maths.warwick.ac.uk/agt/AGTVol\thevolumenumber/agt-\thevolumenumber-\thepapernumber.abs.html}
\immediate\write\gtoutfile{\noexpand\\}
\immediate\write\gtoutfile{}
\ifx\theasciiabstract\relax
\immediate\write\gtoutfile{\theabstract}\else
\immediate\write\gtoutfile{\theasciiabstract}\fi
\immediate\write\gtoutfile{}
\immediate\write\gtoutfile{\noexpand\\}
\immediate\write\gtoutfile{}
\immediate\closeout\gtoutfile}}  %%% end of definition of \makeheadfile

\def\maketitlepage{\makeagttitle\makeheadfile}
\let\makeshorttitle\maketitlepage
\let\maketitle\maketitlepage

\def\ifplaintex{\expandafter\ifx\csname documentclass\endcsname\relax}

\def\gtp{{\mathsurround=0pt\it $\cal G\mskip-2mu$eometry \&\ 
$\cal T\!\!$opology $\cal P\!$ublications}}  % GT publications

\def\recd{{\small Received:\qua\receiveddate\ifx\reviseddate\relax
\else\qquad Revised:\qua\reviseddate\fi\par}} 

\def\lognumber#1{\def\thelognumber{#1}}
\def\volumenumber#1{\def\thevolumenumber{#1}}
\def\volumeyear#1{\def\thevolumeyear{#1}}
\def\papernumber#1{\def\thepapernumber{#1}}
\def\pagenumbers#1#2{\def\startpage{#1}\def\finishpage{#2}}
\def\published#1{\def\publishdate{#1}}

\def\received#1{\def\receiveddate{#1}}
\def\revised#1{\def\reviseddate{#1}}
\def\accepted#1{\def\accepteddate{#1}}

\long\def\asciiabstract#1{\long\def\theasciiabstract{#1}}

\let\\\par\let\thelognumber\relax\let\thevolumenumber\relax
\let\thepapernumber\relax\let\thevolumeyear\relax\let\startpage\relax
\let\finishpage\relax\let\publishdate\relax\let\receiveddate\relax
\let\reviseddate\relax\let\accepteddate\relax\let\theasciititle\relax
\let\theasciiauthors\relax
\let\theasciiabstract\relax

\let\theasciiemail\relax

\ifplaintex
\font\logobig=cmssbx10 scaled 3836
\font\logomed=cmssbx10 scaled 2557
\else
\font\logobig=cmssbx10 scaled 4200
\font\logomed=cmssbx10 scaled 2800
\fi

\long\def\makeagttitle{   %%% start of definition of \makeagttitle
\count0=\startpage
\agt\hfill      %   Journal title (top left) 
\hbox to 45truept{\vbox to 0pt{\vglue -13truept{\logomed A\kern -.37em{\logobig 
T}\kern -.38em G}\vss}\hss}
\break
{\small Volume \thevolumenumber\ (\thevolumeyear)
\startpage--\finishpage\nl
Published: \publishdate}

\vglue .25truein

{\parskip=0pt\leftskip 0pt plus
1fil\def\\{\par\smallskip}{\Large\bf\thetitle}\par\medskip} \vglue
0.05truein

{\parskip=0pt\leftskip 0pt plus 1fil\def\\{\par}{\sc\theauthors}
\par\medskip}%
 
\vglue 0.03truein 

{\small\leftskip 25truept\rightskip 25truept{\bf Abstract}\stdspace\theabstract

{\bf AMS Classification}\stdspace\theprimaryclass
\ifx\thesecondaryclass\relax\else; \thesecondaryclass\fi\par
{\bf Keywords}\stdspace \thekeywords\par}\vglue 7truept

}   %%%% end of definition of \makeagttitle

\ifplaintex
\font\phead=cmsl9 scaled 950
\font\pnum=cmbx10 scaled 913
\font\pfoot=cmsl9 scaled 950
\headline{\vbox to 0pt{\vskip -4.5mm\line{\small\phead\ifnum
\count0=\startpage ISSN 1472-2739 (on-line) 1472-2747 (printed)
\hfill {\pnum\folio}\else\ifodd\count0\def\\{ }% 
\ifx\theshorttitle\relax\thetitle\else\theshorttitle\fi\hfill{\pnum\folio}
\else\def\\{ and }{\pnum\folio}\hfill\ifx\theshortauthors\relax\theauthors
\else\theshortauthors\fi\fi\fi}\vss}}
\footline{\vbox to 0pt{\vglue 0mm\line{\small\pfoot\ifnum\count0=\startpage
\copyright\ \gtp\hfill\else
\agt, Volume \thevolumenumber\ (\thevolumeyear)\hfill\fi}\vss}}
\else
\font=cmsl9 scaled 1050
\font\lnum=cmbx10 
\font=cmsl9 scaled 1050
\makeatletter
\def\@oddhead{{\small\lhead\ifnum\count0=\startpage ISSN 1472-2739 
(on-line) 1472-2747 (printed)\hfill {\lnum\number\count0}\else\ifodd\count0
\def\\{ }\ifx\theshorttitle\relax \thetitle \else\theshorttitle\fi\hfill
{\lnum\number\count0}\else\def\\{ and }{\lnum\number\count0}
\hfill\ifx\theshortauthors\relax 
\theauthors\else\theshortauthors\fi\fi\fi}}\def\@evenhead{\@oddhead}
\def\@oddfoot{\small\lfoot\ifnum\count0=\startpage\copyright\ \gtp\hfill\else
\agt, Volume \thevolumenumber\ (\thevolumeyear)\hfill\fi}
\def\@evenfoot{\@oddfoot}
\makeatother
\fi
\let\maketitlepage\makeagttitle
\let\makeshorttitle\maketitlepage
\let\maketitle\maketitlepage

\newwrite\gtoutfile
\long\gdef\makeheadfile{  %%% start of definition of \makeheadfile
{\def\\{, }\def\s{ }
\immediate\openout\gtoutfile head.xxx
\immediate\write\gtoutfile{To: math@arxiv.org}
\immediate\write\gtoutfile{Subject: put OR rep NNNNN:ppppp}
\immediate\write\gtoutfile{--text follows this line--}
\immediate\write\gtoutfile{Proxy-for: \ifx\theasciiauthors\relax
\theauthors\else\theasciiauthors\fi\s<\ifx\theasciiemail\relax\theemail\else\theasciiemail\fi>}
\immediate\write\gtoutfile{\noexpand\\}
\immediate\write\gtoutfile{Authors: \ifx\theasciiauthors\relax
\theauthors\else\theasciiauthors\fi}
{\def\\{ }\immediate\write\gtoutfile{Title: \ifx\theasciititle\relax
\thetitle\else\theasciititle\fi}}
\immediate\write\gtoutfile{Subj-class: GT or SG, GR etc}
\immediate\write\gtoutfile{MSC-class: \theprimaryclass\ifx\thesecondaryclass\relax\else, \thesecondaryclass\fi}
\immediate\write\gtoutfile{Journal-ref: Algebr. Geom. Topol. \thevolumenumber\s
(\thevolumeyear) \startpage-\finishpage}
\immediate\write\gtoutfile{Comments: Published by Algebraic and
Geometric Topology at}
\immediate\write\gtoutfile{\s\s\s  http://www.maths.warwick.ac.uk/agt/AGTVol\thevolumenumber/agt-\thevolumenumber-\thepapernumber.abs.html}
\immediate\write\gtoutfile{\noexpand\\}
\immediate\write\gtoutfile{}
\ifx\theasciiabstract\relax
\immediate\write\gtoutfile{\theabstract}\else
\immediate\write\gtoutfile{\theasciiabstract}\fi
\immediate\write\gtoutfile{}
\immediate\write\gtoutfile{\noexpand\\}
\immediate\write\gtoutfile{}
\immediate\closeout\gtoutfile}}  %%% end of definition of \makeheadfile

\def\maketitlepage{\makeagttitle\makeheadfile}
\let\makeshorttitle\maketitlepage
\let\maketitle\maketitlepage

\def\ifplaintex{\expandafter\ifx\csname documentclass\endcsname\relax}

\def\gtp{{\mathsurround=0pt\it $\cal G\mskip-2mu$eometry \&\ 
$\cal T\!\!$opology $\cal P\!$ublications}}  % GT publications

\def\recd{{\small Received:\qua\receiveddate\ifx\reviseddate\relax
\else\qquad Revised:\qua\reviseddate\fi\par}} 

\def\lognumber#1{\def\thelognumber{#1}}
\def\volumenumber#1{\def\thevolumenumber{#1}}
\def\volumeyear#1{\def\thevolumeyear{#1}}
\def\papernumber#1{\def\thepapernumber{#1}}
\def\pagenumbers#1#2{\def\startpage{#1}\def\finishpage{#2}}
\def\published#1{\def\publishdate{#1}}

\def\received#1{\def\receiveddate{#1}}
\def\revised#1{\def\reviseddate{#1}}
\def\accepted#1{\def\accepteddate{#1}}

\long\def\asciiabstract#1{\long\def\theasciiabstract{#1}}

\let\\\par\let\thelognumber\relax\let\thevolumenumber\relax
\let\thepapernumber\relax\let\thevolumeyear\relax\let\startpage\relax
\let\finishpage\relax\let\publishdate\relax\let\receiveddate\relax
\let\reviseddate\relax\let\accepteddate\relax\let\theasciititle\relax
\let\theasciiauthors\relax
\let\theasciiabstract\relax

\let\theasciiemail\relax

\ifplaintex
\font\logobig=cmssbx10 scaled 3836
\font\logomed=cmssbx10 scaled 2557
\else
\font\logobig=cmssbx10 scaled 4200
\font\logomed=cmssbx10 scaled 2800
\fi

\long\def\makeagttitle{   %%% start of definition of \makeagttitle
\count0=\startpage
\agt\hfill      %   Journal title (top left) 
\hbox to 45truept{\vbox to 0pt{\vglue -13truept{\logomed A\kern -.37em{\logobig 
T}\kern -.38em G}\vss}\hss}
\break
{\small Volume \thevolumenumber\ (\thevolumeyear)
\startpage--\finishpage\nl
Published: \publishdate}

\vglue .25truein

{\parskip=0pt\leftskip 0pt plus
1fil\def\\{\par\smallskip}{\Large\bf\thetitle}\par\medskip} \vglue
0.05truein

{\parskip=0pt\leftskip 0pt plus 1fil\def\\{\par}{\sc\theauthors}
\par\medskip}%
 
\vglue 0.03truein 

{\small\leftskip 25truept\rightskip 25truept{\bf Abstract}\stdspace\theabstract

{\bf AMS Classification}\stdspace\theprimaryclass
\ifx\thesecondaryclass\relax\else; \thesecondaryclass\fi\par
{\bf Keywords}\stdspace \thekeywords\par}\vglue 7truept

}   %%%% end of definition of \makeagttitle

\ifplaintex
\font\phead=cmsl9 scaled 950
\font\pnum=cmbx10 scaled 913
\font\pfoot=cmsl9 scaled 950
\headline{\vbox to 0pt{\vskip -4.5mm\line{\small\phead\ifnum
\count0=\startpage ISSN 1472-2739 (on-line) 1472-2747 (printed)
\hfill {\pnum\folio}\else\ifodd\count0\def\\{ }% 
\ifx\theshorttitle\relax\thetitle\else\theshorttitle\fi\hfill{\pnum\folio}
\else\def\\{ and }{\pnum\folio}\hfill\ifx\theshortauthors\relax\theauthors
\else\theshortauthors\fi\fi\fi}\vss}}
\footline{\vbox to 0pt{\vglue 0mm\line{\small\pfoot\ifnum\count0=\startpage
\copyright\ \gtp\hfill\else
\agt, Volume \thevolumenumber\ (\thevolumeyear)\hfill\fi}\vss}}
\else
\font=cmsl9 scaled 1050
\font\lnum=cmbx10 
\font=cmsl9 scaled 1050
\makeatletter
\def\@oddhead{{\small\lhead\ifnum\count0=\startpage ISSN 1472-2739 
(on-line) 1472-2747 (printed)\hfill {\lnum\number\count0}\else\ifodd\count0
\def\\{ }\ifx\theshorttitle\relax \thetitle \else\theshorttitle\fi\hfill
{\lnum\number\count0}\else\def\\{ and }{\lnum\number\count0}
\hfill\ifx\theshortauthors\relax 
\theauthors\else\theshortauthors\fi\fi\fi}}\def\@evenhead{\@oddhead}
\def\@oddfoot{\small\lfoot\ifnum\count0=\startpage\copyright\ \gtp\hfill\else
\agt, Volume \thevolumenumber\ (\thevolumeyear)\hfill\fi}
\def\@evenfoot{\@oddfoot}
\makeatother
\fi
\let\maketitlepage\makeagttitle
\let\makeshorttitle\maketitlepage
\let\maketitle\maketitlepage

\newwrite\gtoutfile
\long\gdef\makeheadfile{  %%% start of definition of \makeheadfile
{\def\\{, }\def\s{ }
\immediate\openout\gtoutfile head.xxx
\immediate\write\gtoutfile{To: math@arxiv.org}
\immediate\write\gtoutfile{Subject: put OR rep NNNNN:ppppp}
\immediate\write\gtoutfile{--text follows this line--}
\immediate\write\gtoutfile{Proxy-for: \ifx\theasciiauthors\relax
\theauthors\else\theasciiauthors\fi\s<\ifx\theasciiemail\relax\theemail\else\theasciiemail\fi>}
\immediate\write\gtoutfile{\noexpand\\}
\immediate\write\gtoutfile{Authors: \ifx\theasciiauthors\relax
\theauthors\else\theasciiauthors\fi}
{\def\\{ }\immediate\write\gtoutfile{Title: \ifx\theasciititle\relax
\thetitle\else\theasciititle\fi}}
\immediate\write\gtoutfile{Subj-class: GT or SG, GR etc}
\immediate\write\gtoutfile{MSC-class: \theprimaryclass\ifx\thesecondaryclass\relax\else, \thesecondaryclass\fi}
\immediate\write\gtoutfile{Journal-ref: Algebr. Geom. Topol. \thevolumenumber\s
(\thevolumeyear) \startpage-\finishpage}
\immediate\write\gtoutfile{Comments: Published by Algebraic and
Geometric Topology at}
\immediate\write\gtoutfile{\s\s\s  http://www.maths.warwick.ac.uk/agt/AGTVol\thevolumenumber/agt-\thevolumenumber-\thepapernumber.abs.html}
\immediate\write\gtoutfile{\noexpand\\}
\immediate\write\gtoutfile{}
\ifx\theasciiabstract\relax
\immediate\write\gtoutfile{\theabstract}\else
\immediate\write\gtoutfile{\theasciiabstract}\fi
\immediate\write\gtoutfile{}
\immediate\write\gtoutfile{\noexpand\\}
\immediate\write\gtoutfile{}
\immediate\closeout\gtoutfile}}  %%% end of definition of \makeheadfile

\def\maketitlepage{\makeagttitle\makeheadfile}
\let\makeshorttitle\maketitlepage
\let\maketitle\maketitlepage

\def\ifplaintex{\expandafter\ifx\csname documentclass\endcsname\relax}

\def\gtp{{\mathsurround=0pt\it $\cal G\mskip-2mu$eometry \&\ 
$\cal T\!\!$opology $\cal P\!$ublications}}  % GT publications

\def\recd{{\small Received:\qua\receiveddate\ifx\reviseddate\relax
\else\qquad Revised:\qua\reviseddate\fi\par}} 

\def\lognumber#1{\def\thelognumber{#1}}
\def\volumenumber#1{\def\thevolumenumber{#1}}
\def\volumeyear#1{\def\thevolumeyear{#1}}
\def\papernumber#1{\def\thepapernumber{#1}}
\def\pagenumbers#1#2{\def\startpage{#1}\def\finishpage{#2}}
\def\published#1{\def\publishdate{#1}}

\def\received#1{\def\receiveddate{#1}}
\def\revised#1{\def\reviseddate{#1}}
\def\accepted#1{\def\accepteddate{#1}}

\long\def\asciiabstract#1{\long\def\theasciiabstract{#1}}

\let\\\par\let\thelognumber\relax\let\thevolumenumber\relax
\let\thepapernumber\relax\let\thevolumeyear\relax\let\startpage\relax
\let\finishpage\relax\let\publishdate\relax\let\receiveddate\relax
\let\reviseddate\relax\let\accepteddate\relax\let\theasciititle\relax
\let\theasciiauthors\relax
\let\theasciiabstract\relax

\let\theasciiemail\relax

\ifplaintex
\font\logobig=cmssbx10 scaled 3836
\font\logomed=cmssbx10 scaled 2557
\else
\font\logobig=cmssbx10 scaled 4200
\font\logomed=cmssbx10 scaled 2800
\fi

\long\def\makeagttitle{   %%% start of definition of \makeagttitle
\count0=\startpage
\agt\hfill      %   Journal title (top left) 
\hbox to 45truept{\vbox to 0pt{\vglue -13truept{\logomed A\kern -.37em{\logobig 
T}\kern -.38em G}\vss}\hss}
\break
{\small Volume \thevolumenumber\ (\thevolumeyear)
\startpage--\finishpage\nl
Published: \publishdate}

\vglue .25truein

{\parskip=0pt\leftskip 0pt plus
1fil\def\\{\par\smallskip}{\Large\bf\thetitle}\par\medskip} \vglue
0.05truein

{\parskip=0pt\leftskip 0pt plus 1fil\def\\{\par}{\sc\theauthors}
\par\medskip}%
 
\vglue 0.03truein 

{\small\leftskip 25truept\rightskip 25truept{\bf Abstract}\stdspace\theabstract

{\bf AMS Classification}\stdspace\theprimaryclass
\ifx\thesecondaryclass\relax\else; \thesecondaryclass\fi\par
{\bf Keywords}\stdspace \thekeywords\par}\vglue 7truept

}   %%%% end of definition of \makeagttitle

\ifplaintex
\font\phead=cmsl9 scaled 950
\font\pnum=cmbx10 scaled 913
\font\pfoot=cmsl9 scaled 950
\headline{\vbox to 0pt{\vskip -4.5mm\line{\small\phead\ifnum
\count0=\startpage ISSN 1472-2739 (on-line) 1472-2747 (printed)
\hfill {\pnum\folio}\else\ifodd\count0\def\\{ }% 
\ifx\theshorttitle\relax\thetitle\else\theshorttitle\fi\hfill{\pnum\folio}
\else\def\\{ and }{\pnum\folio}\hfill\ifx\theshortauthors\relax\theauthors
\else\theshortauthors\fi\fi\fi}\vss}}
\footline{\vbox to 0pt{\vglue 0mm\line{\small\pfoot\ifnum\count0=\startpage
\copyright\ \gtp\hfill\else
\agt, Volume \thevolumenumber\ (\thevolumeyear)\hfill\fi}\vss}}
\else
\font=cmsl9 scaled 1050
\font\lnum=cmbx10 
\font=cmsl9 scaled 1050
\makeatletter
\def\@oddhead{{\small\lhead\ifnum\count0=\startpage ISSN 1472-2739 
(on-line) 1472-2747 (printed)\hfill {\lnum\number\count0}\else\ifodd\count0
\def\\{ }\ifx\theshorttitle\relax \thetitle \else\theshorttitle\fi\hfill
{\lnum\number\count0}\else\def\\{ and }{\lnum\number\count0}
\hfill\ifx\theshortauthors\relax 
\theauthors\else\theshortauthors\fi\fi\fi}}\def\@evenhead{\@oddhead}
\def\@oddfoot{\small\lfoot\ifnum\count0=\startpage\copyright\ \gtp\hfill\else
\agt, Volume \thevolumenumber\ (\thevolumeyear)\hfill\fi}
\def\@evenfoot{\@oddfoot}
\makeatother
\fi
\let\maketitlepage\makeagttitle
\let\makeshorttitle\maketitlepage
\let\maketitle\maketitlepage

\newwrite\gtoutfile
\long\gdef\makeheadfile{  %%% start of definition of \makeheadfile
{\def\\{, }\def\s{ }
\immediate\openout\gtoutfile head.xxx
\immediate\write\gtoutfile{To: math@arxiv.org}
\immediate\write\gtoutfile{Subject: put OR rep NNNNN:ppppp}
\immediate\write\gtoutfile{--text follows this line--}
\immediate\write\gtoutfile{Proxy-for: \ifx\theasciiauthors\relax
\theauthors\else\theasciiauthors\fi\s<\ifx\theasciiemail\relax\theemail\else\theasciiemail\fi>}
\immediate\write\gtoutfile{\noexpand\\}
\immediate\write\gtoutfile{Authors: \ifx\theasciiauthors\relax
\theauthors\else\theasciiauthors\fi}
{\def\\{ }\immediate\write\gtoutfile{Title: \ifx\theasciititle\relax
\thetitle\else\theasciititle\fi}}
\immediate\write\gtoutfile{Subj-class: GT or SG, GR etc}
\immediate\write\gtoutfile{MSC-class: \theprimaryclass\ifx\thesecondaryclass\relax\else, \thesecondaryclass\fi}
\immediate\write\gtoutfile{Journal-ref: Algebr. Geom. Topol. \thevolumenumber\s
(\thevolumeyear) \startpage-\finishpage}
\immediate\write\gtoutfile{Comments: Published by Algebraic and
Geometric Topology at}
\immediate\write\gtoutfile{\s\s\s  http://www.maths.warwick.ac.uk/agt/AGTVol\thevolumenumber/agt-\thevolumenumber-\thepapernumber.abs.html}
\immediate\write\gtoutfile{\noexpand\\}
\immediate\write\gtoutfile{}
\ifx\theasciiabstract\relax
\immediate\write\gtoutfile{\theabstract}\else
\immediate\write\gtoutfile{\theasciiabstract}\fi
\immediate\write\gtoutfile{}
\immediate\write\gtoutfile{\noexpand\\}
\immediate\write\gtoutfile{}
\immediate\closeout\gtoutfile}}  %%% end of definition of \makeheadfile

\def\maketitlepage{\makeagttitle\makeheadfile}
\let\makeshorttitle\maketitlepage
\let\maketitle\maketitlepage

\lognumber{2}
\volumenumber{2}
\volumeyear{2002}
\papernumber{2}
\published{15 January 2002}
\pagenumbers{11}{32}
\received{25 June 2001}
\revised{11 January 2002}
\accepted{11 January 2002}

\usepackage{graphicx}
\usepackage{latexsym}
\usepackage{amscd}
\usepackage{amsmath,amssymb}
\usepackage{psfrag}
\usepackage{hhline}

\newcommand{\ints}{\mathbb{Z}}

\newcommand{\ann}{\mathcal{C}}
\newcommand{\annu}{\mathcal{C}^{(n,p)}}

\newcommand{\del}{\frac{v^{-1}-v}{s-s^{-1}}}
\theoremstyle{plain}
\newtheorem{thm}{Theorem}
\newtheorem{lem}[thm]{Lemma}
\newtheorem{prop}[thm]{Proposition}
\newtheorem*{cor}{Corollary}
\theoremstyle{definition}
\newtheorem{defn}{Definition}
\newtheorem*{exmp}{Example}
\theoremstyle{remark}
\newtheorem*{nota}{Notation}
\newtheorem*{rem}{Remark}
\newtheorem*{obs}{Observation}

\begin{document}

% Header
\title{Homfly Polynomials of Generalized Hopf Links}
\authors{Hugh R. Morton\\Richard J. Hadji}
\address{Department of Mathematical Sciences, University of Liverpool\\Peach Street, Liverpool, L69 3ZL, UK}
\email{morton@liv.ac.uk, rhadji@liv.ac.uk}
\url{http://www.liv.ac.uk/\char'176su14/knotgroup.html}

% Abstract
\begin{abstract}
Following the recent work by T.-H. Chan in \cite{Chan00} on reverse string parallels of the Hopf link we give an alternative approach to finding the Homfly polynomials of these links, based on the Homfly skein of the annulus. We establish that two natural skein maps have distinct eigenvalues, answering a question raised by Chan, and use this result to calculate the Homfly polynomial of some more general reverse string satellites of the Hopf link.
\end{abstract}
\asciiabstract{
Following the recent work by T.-H. Chan in [HOMFLY polynomial of some
generalized Hopf links, J. Knot Theory Ramif. 9 (2000) 865--883] on
reverse string parallels of the Hopf link we give an alternative
approach to finding the Homfly polynomials of these links, based on
the Homfly skein of the annulus. We establish that two natural skein
maps have distinct eigenvalues, answering a question raised by Chan,
and use this result to calculate the Homfly polynomial of some more
general reverse string satellites of the Hopf link.}
% Misc.
\maketitle

% Main Text
\section*{Introduction} In \cite{Chan00} T.-H. Chan discusses the Homfly
polynomial of reverse string parallels $H(k_1,k_2;n_1,n_2)$ of the Hopf link.  In
this paper we analyse their structure more closely using the Homfly skein of the
annulus and identify the eigenvalues and eigenvectors which occur naturally in
this approach.

This allows us to readily calculate the Homfly polynomial of satellites of the
Hopf link which consist of a reverse string parallel around one component
combined with a completely general reverse string decoration on the other.
\subsection*{The Homfly polynomial} Various versions of the Homfly polynomial
appear in the literature.  The framed version to the fore in this paper is
defined by the following skein relations:
\begin{eqnarray*}
{\raisebox{-14pt}{\includegraphics[height=30pt]{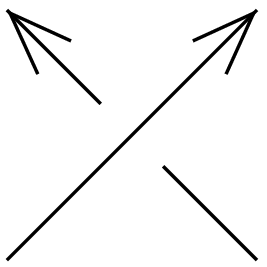}}}-{\raisebox{-14pt}{\includegraphics[height=30pt]{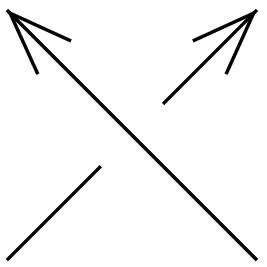}}}&=&(s-s^{-1}){\raisebox{-14pt}{\includegraphics[height=30pt]{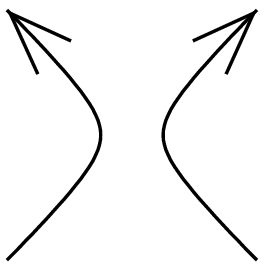}}},
\\
\text{and}\phantom{--}{\raisebox{-14pt}{\includegraphics[height=30pt]{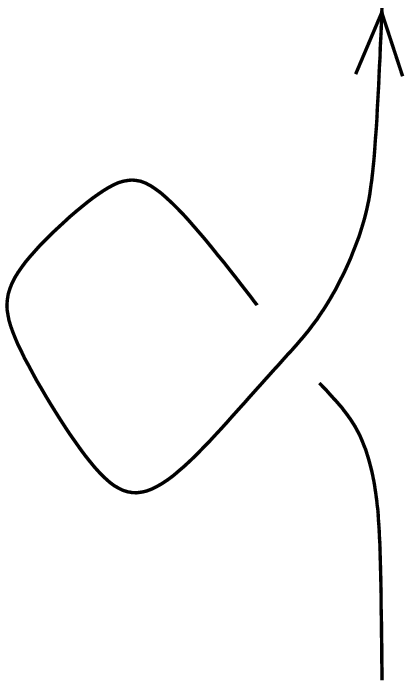}}}&=&v^{-1}\,{\raisebox{-14pt}{\includegraphics[height=30pt]{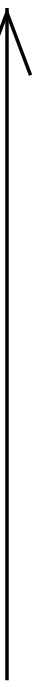}}}\,
,
\end{eqnarray*} with the Homfly polynomial of the empty diagram being normalized
to 1, and that of the null-homotopic loop therefore being $\delta=\del$.  There
is a discussion of isomorphic variants of these skein relations given in
\cite{AM98} and \cite{Morton01}.
\begin{nota} For a link $L$, we denote the evaluation of its Homfly polynomial by
$P(L)$.
\end{nota}
\begin{rem} (i)\qua The Homfly polynomial of the $m$-component unlink,
$\mathcal{U}^m=$\break $\sqcup_{i=1}^m
{\raisebox{-4pt}{\includegraphics[height=15pt]{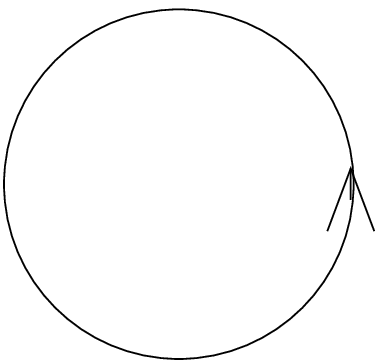}}}\,$, is
$P(\mathcal{U}^m)=\delta^m$.

$\phantom{Remar.....}$(ii)\qua If $L^*$ is the reflection of a link $L$, then 
\[ P(L^*)(s,v)=P(L)(s^{-1},v^{-1}).\]
\end{rem}
\subsection*{Diagrams in the annulus} We shall now introduce the basic idea of
\emph{diagrams in the annulus}.  Given the annulus $F=S^1\times I$, a diagram in
$F$ consists of closed curves (as with a standard knot diagram) with a finite
number of crossing points.  At a crossing point the strands are distinguished in
the conventional way as an over-crossing and an under-crossing.
\subsection*{Satellites of Hopf links} 
\begin{figure}[ht!]
\psfrag{1}{$H_+$}
\psfrag{2}{$H_-$}
\begin{center}
\includegraphics[width=0.5\textwidth]{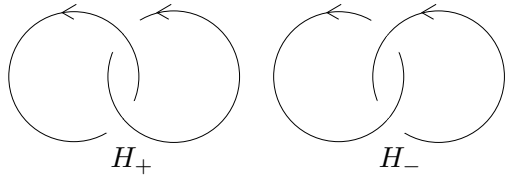}
\caption{The links $H_+$ and $H_-$}
\label{fig:hopf}
\end{center}
\end{figure}

The Hopf link is the simplest non-trivial
link involving just two unknots linked together.  When giving this link
orientation, two distinct links are formed.  We shall call these $H_+$ and $H_-$,
as shown in Figure~\ref{fig:hopf}.
 The Homfly polynomial can then be calculated using the above skein
relations.  We have that:
\begin{eqnarray*} P(H_+)&=&\left(\frac{v^{-1}-v}{s-s^{-1}}\right)^2+v^{-2}-1; \\
\text{and}\phantom{-}P(H_-)&=&\left(\frac{v^{-1}-v}{s-s^{-1}}\right)^2+v^2-1.
\end{eqnarray*}

\begin{figure}[ht!]
\psfrag{1}{$P_1$}
\psfrag{2}{$P_2$}
\begin{center}
\includegraphics[width=0.25\textwidth]{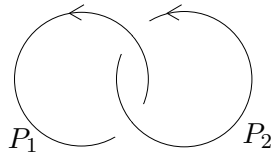}
\caption{The link $H_+(P_1,P_2)$}
\label{fig:decor}
\end{center}
\end{figure}

We now use $H_+$ and $H_-$ as starting points for the construction of satellite
links.  We do this by considering the two components of the Hopf links and
decorating them.  For example, take $P_1$ and $P_2$ as diagrams in the annulus. 
Now starting with $H_+$ we decorate its two components with $P_1$ and $P_2$
respectively, obtaining the link $H_+(P_1,P_2)$, as shown in
Figure~\ref{fig:decor}.  Now clearly $H_+(P_1,P_2)$ and $H_+(P_2,P_1)$ are
equivalent links.  An analogous construction is possible for $H_-$.

\begin{figure}[ht!]
\psfrag{1}{$\scriptstyle{n_2}$}
\psfrag{2}{$\scriptstyle{n_1}$}
\psfrag{3}{$\scriptstyle{k_2}$}
\psfrag{4}{$\scriptstyle{k_1}$}
\begin{center}
\includegraphics[width=0.5\textwidth]{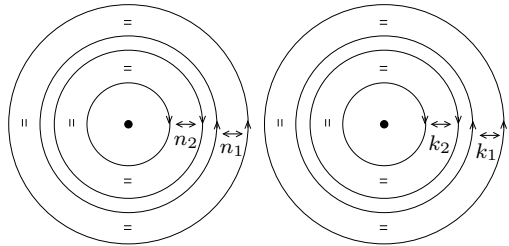}
\caption{The diagrams $P_1$ and $P_2$}
\label{fig:decor2}
\end{center}
\end{figure}

\begin{figure}[ht!]
\psfrag{1}{$n_1$}
\psfrag{2}{$n_2$}
\psfrag{3}{$k_1$}
\psfrag{4}{$k_2$}
\begin{center}
\includegraphics[width=0.45\textwidth,height=6cm]{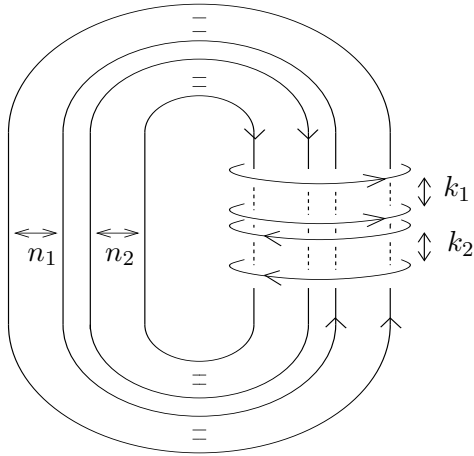}
\caption{The generalized Hopf link $H(k_1,k_2;n_1,n_2)$}
\label{fig:genhopf}
\end{center}
\end{figure}

With such a construction, a great variety of links may be realised.  In
particular, the generalized Hopf links which are the topic of \cite{Chan00} can
be constructed.  For example, if we take $P_1$ and $P_2$ as shown in
Figure~\ref{fig:decor2},
 then $H_+(P_1,P_2)$ is the link Chan refers to as
$H(k_1,k_2;n_1,n_2)$, as shown in Figure~\ref{fig:genhopf}.

With such links in mind, we make the following observation:
\begin{obs} The links \[H(k_1,k_2;n_1,n_2), H(n_1,n_2;k_1,k_2),
H(k_2,k_1;n_2,n_1), H(n_2,n_1;k_2,k_1),\] and \[H^*(k_2,k_1;n_1,n_2),
H^*(n_1,n_2;k_2,k_1), H^*(k_1,k_2;n_2,n_1), H(n_2,n_1;k_1,k_2),\] are all
\emph{equivalent links}.
\end{obs}
\section{The Skein of the Annulus} We have introduced the concept of having
diagrams in the annulus, and used this to construct satellites of the Hopf
links.  We now describe the \emph{skein of the annulus}, denoted by $\ann$.

\begin{figure}[ht!]
\psfrag{1}{$X$}
\begin{center}
\includegraphics[width=0.1\textwidth]{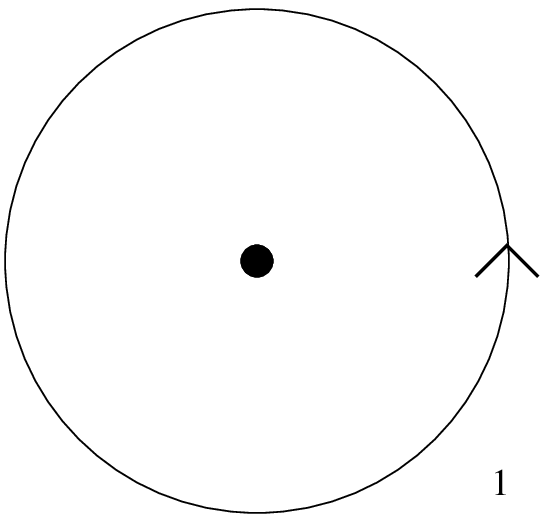}
\caption{An element $X\in\ann$}
\label{fig:X}
\end{center}
\end{figure}

The Homfly skein of the annulus $\ann$, as discussed in \cite{Morton93} and
originally in the preprint of \cite{Turaev97} in 1988, is defined as linear
combinations of diagrams in the annulus, modulo the Homfly skein relations given
above in the Introduction.  We shall represent an element $X\in\ann$
diagramatically as in Figure~\ref{fig:X}.

The skein $\ann$ has a product induced by placing one annulus outside another. 
This defines a bilinear product $\ann\times\ann\rightarrow\ann$, under which
$\ann$ becomes an algebra.  This algebra is clearly commutative (lift the inner
annulus up and stretch it so that the outer one will fit inside it).

\begin{figure}[ht!]
\psfrag{1}{$\scriptstyle{m-1}$}
\begin{center}
\includegraphics[width=0.3\textwidth]{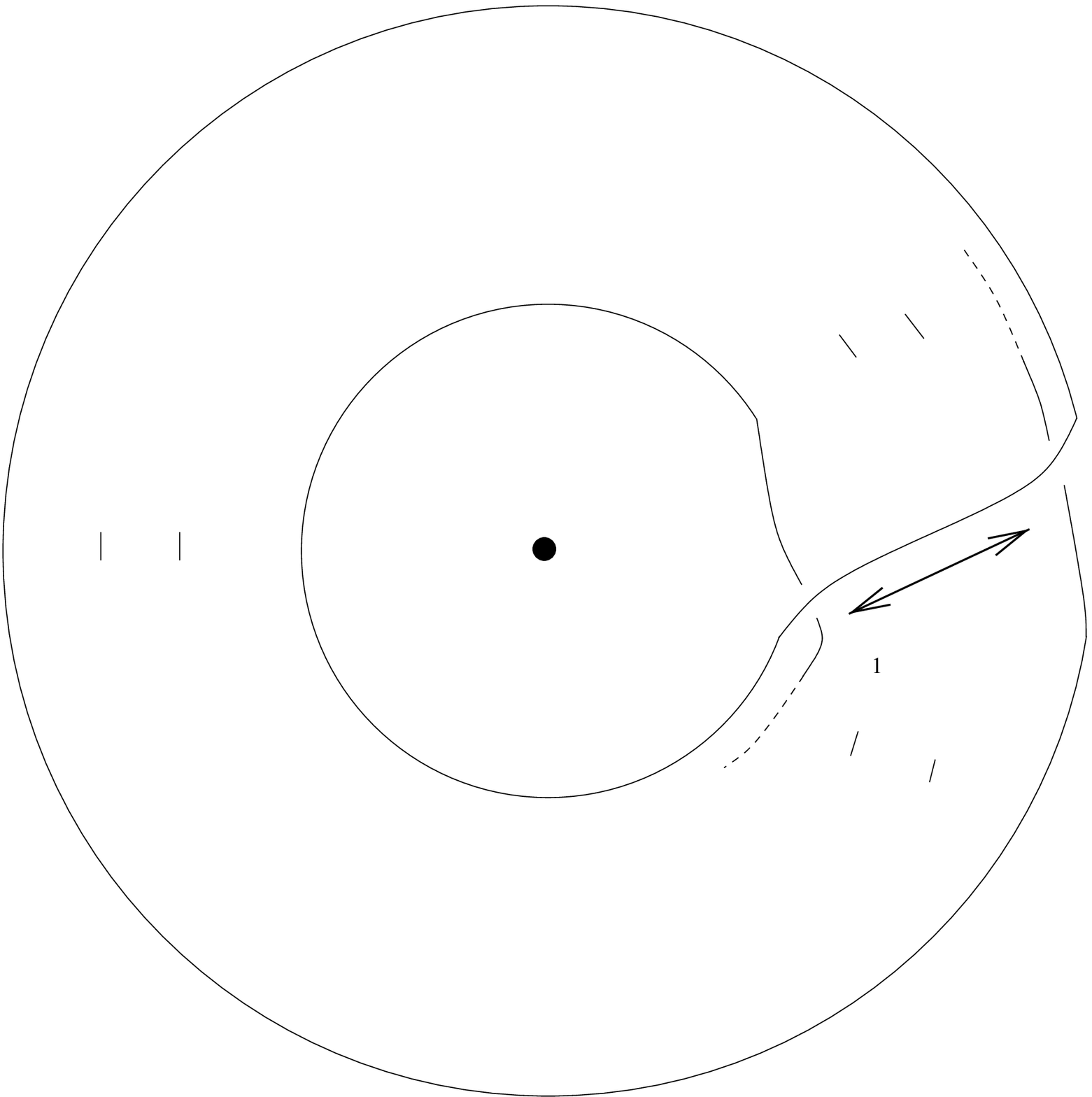}
\caption{An element $A_m\in\ann$, for $m\in\ints$}
\label{fig:Ai}
\end{center}
\end{figure}

Turaev \cite{Turaev97} showed that $\ann$ is freely generated as an algebra by
$\{A_m,m\in\ints\}$ where $A_m$ is the skein element shown in Figure~\ref{fig:Ai}
and the sign of the index $m$ indicates the orientation of the curve.  A positive
$m$ denotes counter-clockwise orientation and a negative $m$ clockwise
orientation.  The element $A_0$ is the identity element, represented by the empty
diagram.

We now define two natural \emph{linear} maps, $\varphi$ and $\bar{\varphi}$, on
the skein of the annulus in the following way:
\begin{eqnarray*}
\varphi:\ann&\rightarrow&\ann \\
{\raisebox{-16pt}{\includegraphics[height=35pt]{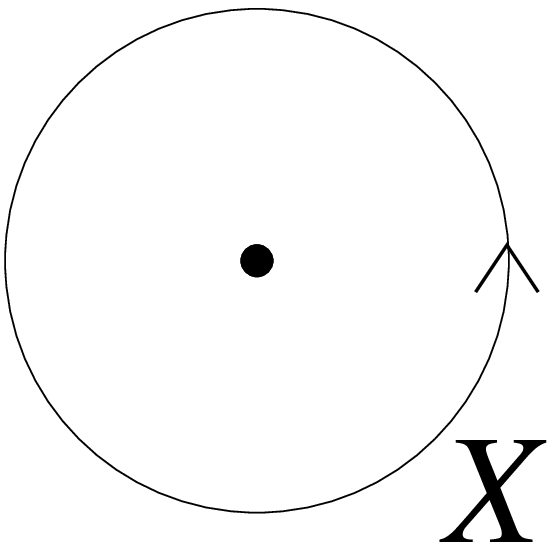}}}&\mapsto& {\raisebox{-16pt}
{\includegraphics[height=35pt]{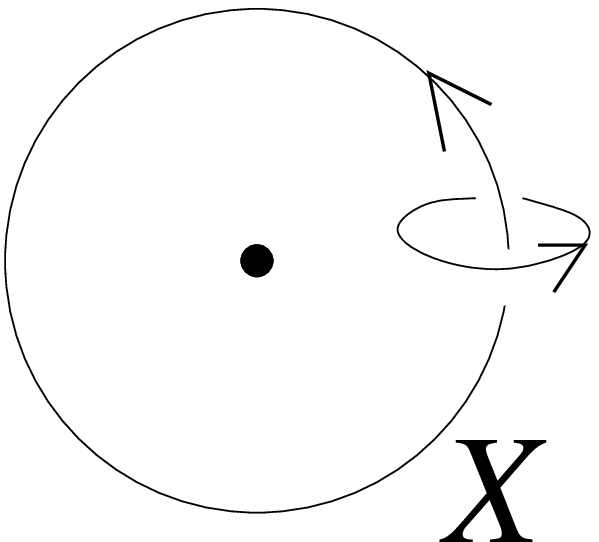}}}, \\
\text{and}\phantom{-}\bar{\varphi}:\ann&\rightarrow&\ann \\
{\raisebox{-16pt}{\includegraphics[height=35pt]{X1}}}&\mapsto& {\raisebox{-16pt}
{\includegraphics[height=35pt]{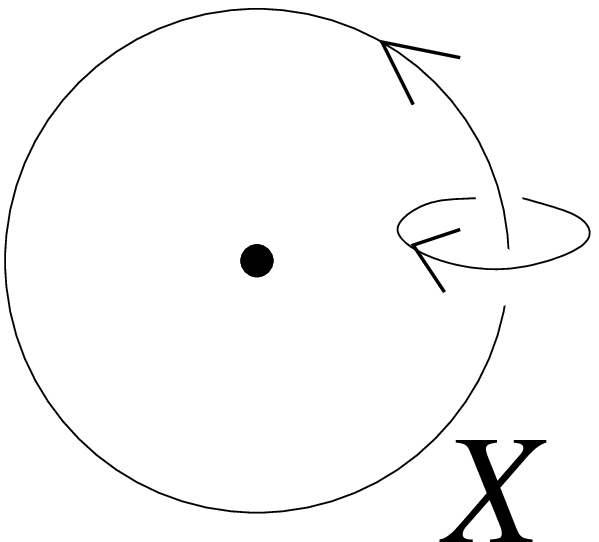}}}.
\end{eqnarray*}

These two maps can be related via a map $\rho:\ann\rightarrow\ann$ which takes
the annulus and its contents and \emph{flips it over}.  Clearly
$\rho^{-1}=\rho$.  It is then clear that
$\bar{\varphi}=\rho^{-1}\varphi\rho=\rho\,\varphi\rho$.  However, this does not
mean that the maps $\varphi$ and $\bar{\varphi}$ are conjugate as this does not
define an inner automorphism of $\ann$.

Now consider the satellites of Hopf links discussed earlier in the Introduction
as elements of the skein $\ann$.  We can then use compositions of the maps
$\varphi$ and $\bar{\varphi}$ to construct a subset of such links.  In
particular, for the element $A=A_1^{n_1}A_{-1}^{n_2}\in\ann$, we have:
\[H(k_1,k_2;n_1,n_2)=\varphi^{k_1}(\bar{\varphi}^{k_2}(A)).\]

It will therefore aid our investigation of the $H(k_1,k_2;n_1,n_2)$ and their
Homfly polynomial if we were to look more closely at the  maps $\varphi$ and
$\bar{\varphi}$, in particular at their eigenvalues.  This will be achieved
through considering certain subspaces of $\ann$ and the restrictions of the maps
$\varphi$ and $\bar{\varphi}$ to these subspaces.

\section{Subspaces of $\ann$} The algebra $\ann$ can be thought of as the product
of subalgebras $\ann^+$ and $\ann^-$ which are generated by $\{A_m:m\in\ints
,m\ge 0\}$ and $\{A_m:m\in\ints ,m\le 0\}$ respectively.

\begin{rem} In his thesis, \cite{Lukac}, Lukac shows how to calculate the Homfly
polynomial of any satellite of the Hopf link, when the decorations are chosen
from $\ann^+$. Here we consider the full skein $\ann$, allowing more general
reverse string decorations.
\end{rem}

\subsection{The subspace $\ann^{(n)}\subset\ann^+$} 
\begin{figure}[ht!]
\psfrag{1}{$T$}
\psfrag{2}{$n$}
\begin{center}
\includegraphics[width=0.12\textwidth]{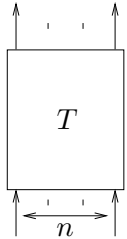}
\caption{An oriented $n$-tangle}
\label{fig:Hn}
\end{center}
\end{figure}

Given linear combinations of
oriented $n$-tangles as shown in Figure~\ref{fig:Hn},
 modulo the Homfly skein relations, we can form an algebra. 
Multiplication within this algebra is induced by stacking one tangle on top of
another.

Now take the well-known Hecke algebra $H_n$ of type $A_{n-1}$.  It has several
different incarnations, but is most conveniently thought of in this context as
having explicit presentation:
\[
H_n=\left<\sigma_i:i=1,\ldots,n-1\left|\begin{array}{l}\sigma_i\sigma_j=\sigma_j\sigma_i:|i-j|>1;\\\sigma_i\sigma_{i+1}\sigma_i=\sigma_{i+1}\sigma_i\sigma_{i+1}:1\le
i<n-1;\\\sigma_i-\sigma_i^{-1}=z\end{array}\right.\right>.
\]

It is shown in \cite{MT90} that $H_n$, with $z=s-s^{-1}$ and coefficient ring
extended to include $v^{\pm 1}$ and $s^{\pm 1}$, is isomorphic to the skein
theoretic algebra described above. 
In this algebra the extra variable $v$ in the coefficient ring allows us to reduce general
tangles to linear combinations of braids, by means of the skein relations in the introduction.
The variable $v$ comes into play in dealing with curls using the second skein relation and in handling disjoint closed curves.

\begin{figure}[ht!]
\psfrag{1}{$T$}
\psfrag{2}{$n$}
\begin{center}
\includegraphics[width=0.25\textwidth]{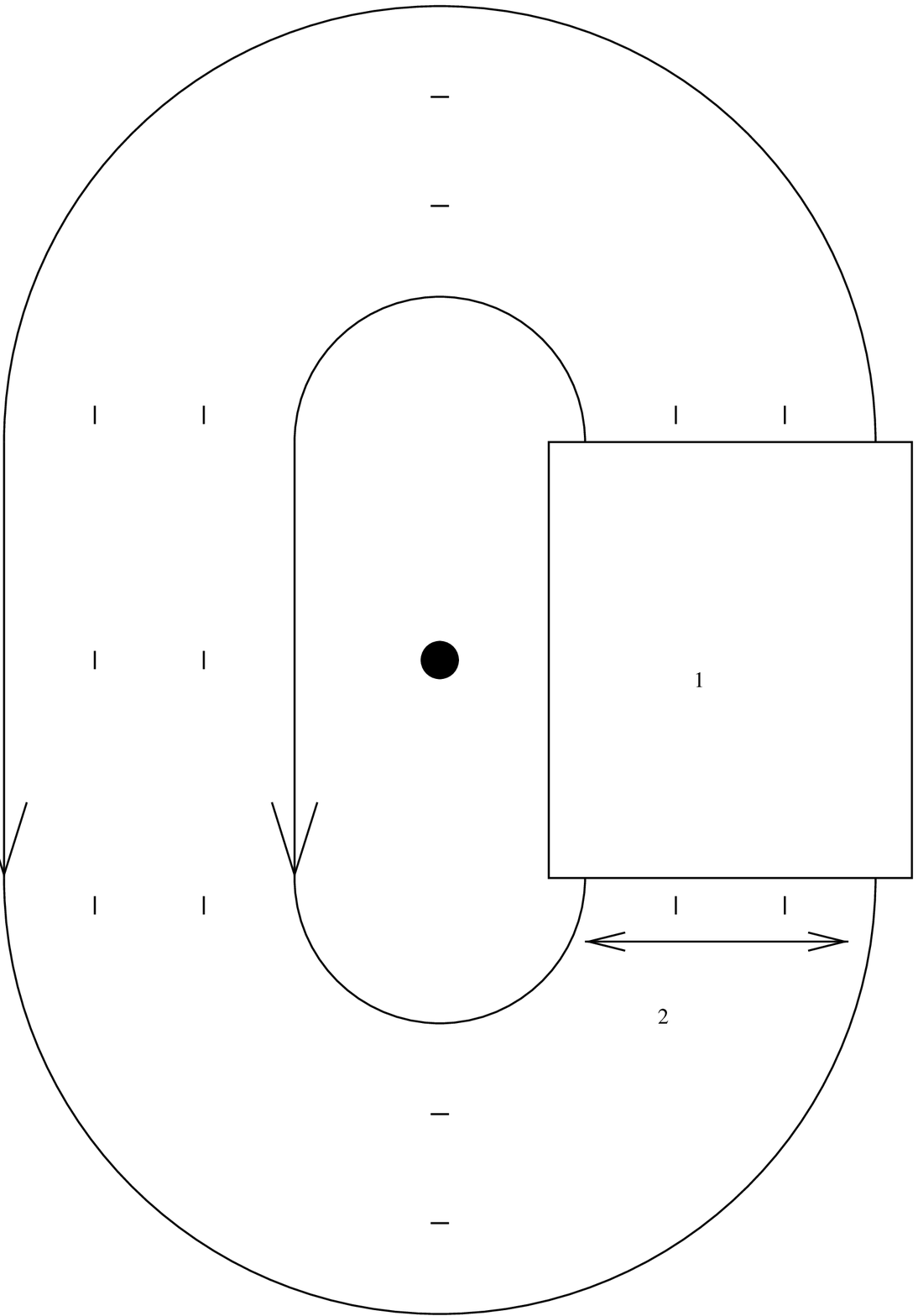}
\caption{An element of the subspace $\ann^{(n)}$}
\label{fig:Cn}
\end{center}
\end{figure}

Wiring these $n$-tangles into the annulus as shown in Figure~\ref{fig:Cn}
 gives a linear subspace of $\ann^+$ which we shall call
$\ann^{(n)}$.  This subspace is the image of $H_n$ under the closure map
$\wedge:H_n\rightarrow\ann^{(n)}$.  For an $n$-tangle $T\in H_n$, we denote its
image under the closure map $\wedge(T)$ or $\hat{T}$.

The subspace $\ann^{(n)}$ is then spanned by monomials in $\{A_m\}$, with
$m\in\ints^+$, of total weight $n$, where $\text{wt}(A_m)=m$.  It is clear that
this spanning set consists of $\pi(n)$ elements, the number of partitions of
$n$.  (The standard notation used for the number of partitions of an integer $n$
is $p(n)$; our alternative has been chosen to avoid a clash with notation
required later in this paper.)  $\ann^+$ is then graded as an algebra:
\[\ann^+=\bigoplus_{n=0}^\infty\ann^{(n)}.\]

Now due to the relationship between $H_n$ and $\ann^{(n)}$ it will be useful here
to recall some well-established facts about $H_n$.  In particular, we shall
concentrate on facts about certain elements in $H_n$.

Firstly, there is a set of quasi-idempotent elements of $H_n$ discussed first in
\cite{Gyoja86} and given a geometric interpretation in \cite{AistonPhD} (see also
\cite{AM98}).  We shall denote these elements $e_\lambda$, one for each partition
$\lambda$ of $n$, with $\emptyset$ denoting the unique partition of $0$.

\begin{figure}[ht!]
\psfrag{1}{$n$}
\psfrag{2}{$H_-$}
\begin{center}
\includegraphics[width=0.2\textwidth]{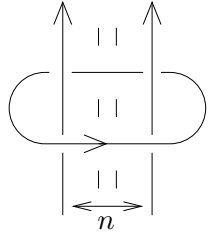}
\caption{The $n$-tangle $T^{(n)}$}
\label{fig:Tn}
\end{center}
\end{figure}

Now, given the element $T^{(n)}\in H_n$ shown in Figure~\ref{fig:Tn}, one can use
skein theoretic techniques to prove the following corollary of Theorem~19 in
\cite{AM98} (see \cite{Morton01}),

\begin{cor}[of Theorem~19, \cite{AM98}]
$T^{(n)}e_\lambda=t_\lambda e_\lambda$ where
\[ t_\lambda=(s-s^{-1})v^{-1}\sum_{\substack{\text{cells}\\ \text{in
}\lambda}}s^{2(\text{content})}+\delta. \] Moreover, the scalars $t_\lambda$ are
different for each partition $\lambda$.
\end{cor}

Reverse the orientation of the encircling string in $T^{(n)}$ and call this
$\bar{T}^{(n)}$.  Then, using similar techniques, one can show
\begin{lem}
$\bar{T}^{(n)}e_\lambda=\bar{t}_\lambda e_\lambda$ where
\[ \bar{t}_\lambda=-(s-s^{-1})v\sum_{\substack{\text{cells}\\ \text{in
}\lambda}}s^{-2(\text{content})}+\delta. \] Moreover, the scalars
$\bar{t}_\lambda$ are also different for each partition $\lambda$.
\end{lem}

\begin{rem} An alternative proof to this lemma could be made using the natural
skein mirror map $\overline{\phantom{;;}}:\ann^{(n)}\rightarrow\ann^{(n)}$.  This
map switches all crossings in a tangle and inverts the scalars $v$ and $s$ in the
coefficient ring.  Clearly this map can be seen to leave the skein relations
unchanged.  One must then note that the $e_\lambda$ are invariant under this map
and that $\overline{\phantom{;;}}(T^{(n)})=\bar{T}^{(n)}$.  We then apply these
facts to the above Corollary and the result follows immediately.
\end{rem}

We now link these facts to the maps $\varphi$ and $\bar{\varphi}$, through use of
the closure map.  Take an element $S\in H_n$ with $\hat{S}\in\ann^{(n)}$ and
compose it with $T^{(n)}$.  Then $\wedge(ST^{(n)})=\varphi(\hat{S})$.  Similarly
$\wedge(S\bar{T}^{(n)})=\bar{\varphi}(\hat{S})$.

The restrictions $\varphi|_{\ann^{(n)}}$ and $\bar{\varphi}|_{\ann^{(n)}}$
clearly carry $\ann^{(n)}$ to itself.
\begin{thm}\label{thm:mort}{\rm \cite{Morton01}}\qua The eigenvalues of
$\varphi|_{\ann^{(n)}}$ are all distinct as are the eigenvalues of
$\bar{\varphi}|_{\ann^{(n)}}$.
\end{thm}
\begin{proof} We prove the first statement with the second following in exactly
the same way.

Set $Q_\lambda=\hat{e}_\lambda\in\ann^{(n)}$.  Then the closure of
$T^{(n)}e_\lambda$ is $\varphi(Q_\lambda)$.  However, $T^{(n)}e_\lambda=t_\lambda
e_\lambda$, hence $\varphi(Q_\lambda)=t_\lambda Q_\lambda$.  The element
$Q_\lambda$ is then an eigenvector of $\varphi$ with eigenvalue $t_\lambda$. 
There are $\pi(n)$ of these eigenvectors, and the eigenvalues are all distinct by
\cite{AM98}.  Since $\ann^{(n)}$ is spanned by $\pi(n)$ elements we can deduce
that the elements $Q_\lambda$ form a basis for $\ann^{(n)}$ and that the
eigenspaces are all $1$-dimensional.
\end{proof} This proof is quite instructive as it establishes that the
$Q_\lambda$ with $|\lambda|=n$ are a basis for $\ann^{(n)}$.  Hence any element
in $\ann^{(n)}$ can be written as a linear combination of the $Q_\lambda$ with
$|\lambda|=n$.  It also follows that \emph{any} element of $\ann^{(n)}$ which is
an eigenvector of $\varphi$ (and similarly $\bar{\varphi}$) must be a multiple of
some $Q_\lambda$.  Finally, we notice that the eigenvalues of the $\varphi$ and
$\bar{\varphi}$ are the $t_\lambda$ and $\bar{t}_\lambda$ we found earlier.

\subsection{The subspace $\ann^{(n,p)}$} 
\begin{figure}[ht!]
\psfrag{1}{$T$}
\psfrag{2}{$n$}
\psfrag{3}{$p$}
\begin{center}
\includegraphics[width=0.17\textwidth]{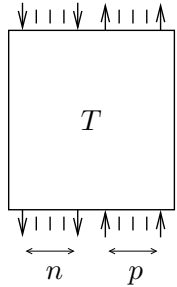}
\caption{An oriented $(n,p)$-tangle}
\label{fig:Mnp}
\end{center}
\end{figure}

We now extend our view of the skein of
the annulus to include strings oriented in both directions.  We do this through
considering the closure of oriented $(n,p)$-tangles such as the one shown in
Figure~\ref{fig:Mnp}.
 We denote the algebra formed through considering linear combinations
of such $(n,p)$-tangles, modulo the Homfly skein relations, by $M_{n,p}$.  For
further information on $M_{n,p}$ see \cite{MW} or \cite{Hadji99}.  The image of
$M_{n,p}$ under the closure map shall be denoted $\ann^{(n,p)}\subset\ann$.

Unlike the case for $\ann^{(n)}$ where $\ann^{(n)}\cap\ann^{(n-1)}=\emptyset$, we
have that:
\[
\annu\supset\ann^{(n-1,p-1)}\supset\ann^{(n-2,p-2)}\supset\cdots\supset\left\{\begin{array}{ll}
\ann^{(n-p,0)} & \text{if }\min(n,p)=p, \\
\ann^{(0,p-n)} & \text{if }\min(n,p)=n, \end{array}\right.
\] however, it should be noted that for each $\ann^{(i,j)}$ in the sequence, the
difference $i-j$ remains constant.  Also,
\begin{eqnarray*}
\ann^{(m,0)}&\cong &\ann^{(m)}_{(-)} \\
\text{and}\phantom{--}\ann^{(0,m)}&\cong &\ann^{(m)}_{(+)},
\end{eqnarray*} where the subscripts indicate the direction of the strings around
the centre of the annulus.  However, we do have that
$\ann^{(n_1,p_1)}\cap\ann^{(n_2,p_2)}=\emptyset$ if $n_1-p_1\ne n_2-p_2$.

We find that $\ann^{(n,p)}$ is spanned by suitably weighted monomials in
\[\{A_{-n},\ldots,A_{-1},A_0,A_1,\ldots,A_p\}.\]  We can then see that:
\[
 \annu=\left(\ann^{(n)}_{(-)}\times\ann^{(p)}_{(+)}\right)+\ann^{(n-1,p-1)}.
\]

The spanning set of $\ann^{(n,p)}$ consists of $\pi(n,p)$ elements, where
\begin{eqnarray*}
\pi(n,p)&:=&\sum_{j=0}^{k}\pi(n-j)\pi(p-j)\\
        &=&\pi(n)\pi(p)+\cdots+\pi(n-k)\pi(p-k),
\end{eqnarray*} with $k=\min(n,p)$.

Similar to the grading of $\ann^+$ with the $\ann^{(n)}$ we can think of the
whole of $\ann$ in terms of the $\ann^{(n,p)}$:
\[
\ann=\bigoplus_{k=-\infty}^{\infty}\left(\bigcup_{n,p\ge 0}\{\annu:n-p=k\}\right).
\]

We now use an example to illustrate what we meant by ``suitably weighted''
monomials in the $A_i$.
\begin{exmp} Consider when $n=3$ and $p=2$.  The spanning set of $\ann^{(3,2)}$
consists of $9\,(=3\cdot 2+2\cdot 1+1\cdot 1)$ elements, since
\[
\ann^{(3,2)}=\left(\ann^{(3)}_{(-)}\times\ann^{(2)}_{(+)}\right)+\left(\ann^{(2)}_{(-)}\times\ann^{(1)}_{(+)}\right)+\left(\ann^{(1)}_{(-)}\times\ann^{(0)}_{(+)}\right).
\] The spanning set is therefore:
\begin{eqnarray*}
\{A_{-3}A_2,\,A_{-3}A_1^2,\,A_{-2}A_{-1}A_2,\,A_{-2}A_{-1}A_1^2,\,A_{-1}^3A_2,\,A_{-1}^3A_1^2,\,A_{-2}A_1,\,\,\,\,\,\,&
& \\ A_{-1}^2A_1,\,A_{-1}\}, & &
\end{eqnarray*} where, for example, the element $A_{-2}A_1$ is obtained from
closing an element in $M_{3,2}$ as shown in Figure~\ref{fig:ex1}.
\begin{figure}[ht!]
\psfrag{1}{$=$}
\psfrag{2}{$\delta$}
\begin{center}
\includegraphics[width=0.5\textwidth]{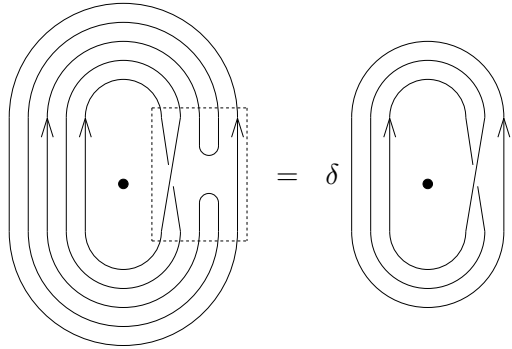}
\caption{The generator $A_{-2}A_1$}
\label{fig:ex1}
\end{center}
\end{figure}
\end{exmp}

\begin{figure}[ht!]
\psfrag{1}{$n$}
\psfrag{2}{$p$}
\begin{center}
\includegraphics[width=0.3\textwidth]{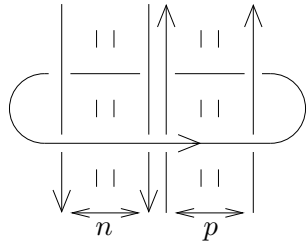}
\caption{The $(n,p)$-tangle $T^{(n,p)}$}
\label{fig:Tnp}
\end{center}
\end{figure}

\begin{figure}[ht!]
\psfrag{1}{$n$}
\psfrag{2}{$p$}
\psfrag{3}{$T_1$}
\psfrag{4}{$\scriptstyle{n-l}$}
\psfrag{5}{$\scriptstyle{p-l}$}
\psfrag{6}{$T_2$}
\begin{center}
\includegraphics[width=0.25\textwidth]{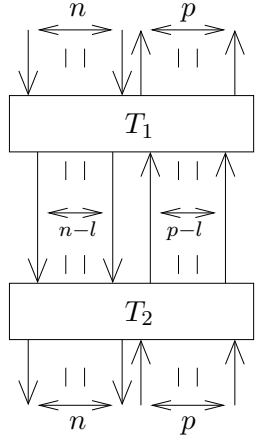}
\caption{A tangle with at least $l$ pairs of strings which \emph{turn back}}
\label{fig:turnback}
\end{center}
\end{figure}

Following an exactly analogous procedure in $M_{n,p}$ as in $H_n$, we define
central elements $T^{(n,p)}$ and $\bar{T}^{(n,p)}$ by diagrams similar to
$T^{(n)}$ and $\bar{T}^{(n)}$ respectively as in Figure~\ref{fig:Tnp}.
\begin{defn}(see \cite{MW},\cite{Hadji99}) Let $M_{n,p}^{(i)}$ denote the sub-algebra of $M_{n,p}$ spanned by elements with ``at least'' $i$ pairs of strings turning back.
\end{defn}
\begin{rem} (i) An $(n,p)$-tangle is said to have ``at least'' $l$ pairs of
strings which \emph{turn back} if it can be written as a product $T_1T_2$ of an
$\{(n,p),(n-l,p-l)\}$-tangle $T_1$ and an $\{(n-l,p-l),(n,p)\}$-tangle $T_2$ as
illustrated in Figure~\ref{fig:turnback}.

$\phantom{Remar....}$(ii) The $M_{n,p}^{(i)}$ are two-sided ideals and there is a
filtration:
\[ M_{n,p}\cong M_{n,p}^{(0)}\rhd M_{n,p}^{(1)}\rhd\cdots\rhd M_{n,p}^{(k)},
\] where $k=\text{min}(n,p)$.
\end{rem}

\begin{lem}\label{thm:decomp}
\begin{eqnarray*} T^{(n,p)}&=&T^{(n,p)\prime}+w,\\
\bar{T}^{(n,p)}&=&\bar{T}^{(n,p)\prime}+\bar{w},
\end{eqnarray*} where
\begin{eqnarray*} T^{(n,p)\prime}&=&T^{(n)}_{(-)}\otimes
1^{(p)}_{(+)}+1^{(n)}_{(-)}\otimes T^{(p)}_{(+)}-\delta\phantom{.}
1^{(n)}_{(-)}\otimes 1^{(p)}_{(+)},\\
\bar{T}^{(n,p)\prime}&=&\bar{T}^{(n)}_{(-)}\otimes
1^{(p)}_{(+)}+1^{(n)}_{(-)}\otimes \bar{T}^{(p)}_{(+)}-\delta\phantom{.}
1^{(n)}_{(-)}\otimes 1^{(p)}_{(+)},
\end{eqnarray*} and $w,\bar{w}\in M_{n,p}^{(1)}$.
\end{lem}
\begin{nota} The tensor product $S\otimes T$ indicates the juxtaposition of
tangles $S$ and $T$.
\end{nota}
\begin{proof}(of Lemma~\ref{thm:decomp}) We prove the result for $T^{(n,p)}$, with the result for
$\bar{T}^{(n,p)}$ following in exactly the same way.  Throughout this proof, we
use a standard notation setting $s-s^{-1}=z$.

We first define some elements in $M_{n,p}$ represented by tangles as shown in
Figure~\ref{fig:TjAj}.
\begin{figure}[ht!]
\psfrag{1}{$1$}
\psfrag{2}{$j$}
\psfrag{3}{$p$}
\psfrag{4}{$1$}
\psfrag{5}{$j$}
\psfrag{6}{$p$}
\psfrag{7}{$\cdots$}
\psfrag{8}{$A(j):=$}
\psfrag{9}{$T(j):=$}
\begin{center}
\includegraphics[width=\textwidth]{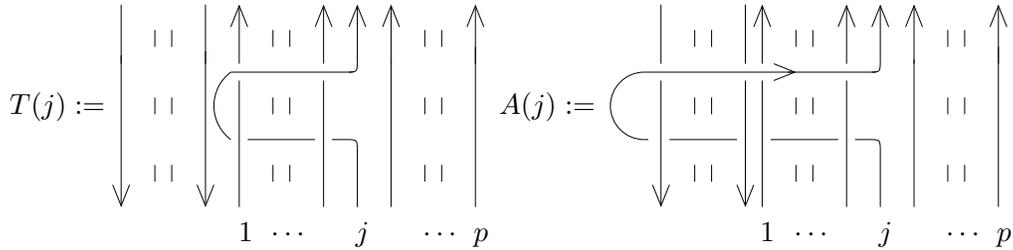}
\caption{The tangles representing the elements $T(j)$ and $A(j)$ for $1\le j\le
p$}
\label{fig:TjAj}
\end{center}
\end{figure}

Now applying the skein relation once to $T^{(n,p)}$ we obtain:
\begin{eqnarray*}
{\raisebox{-16pt}{\includegraphics[height=35pt]{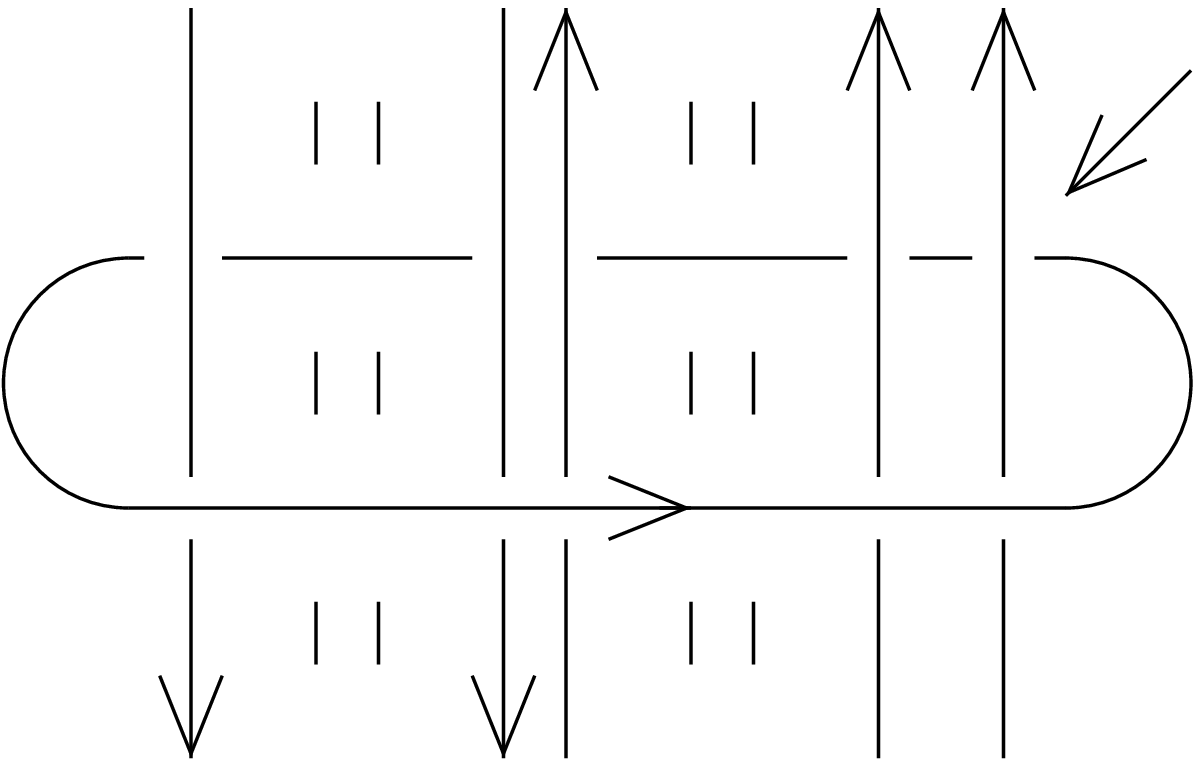}}}&=&{\raisebox{-16pt}{\includegraphics[height=35pt]{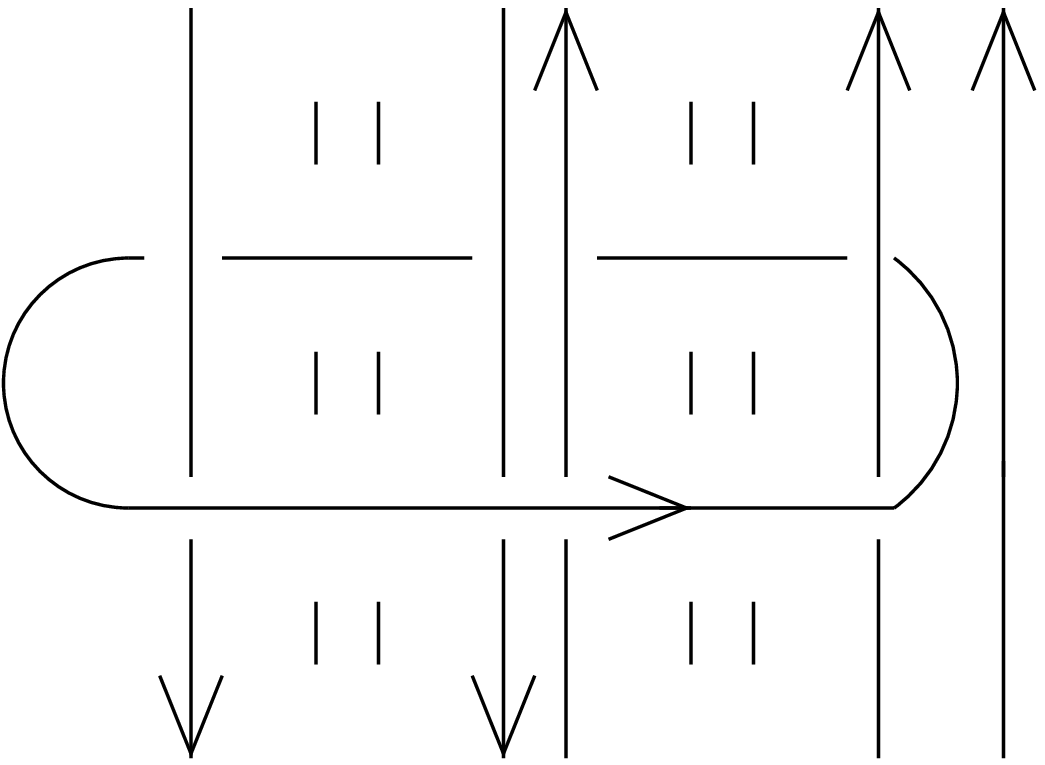}}}+z\,{\raisebox{-16pt}{\includegraphics[height=35pt]{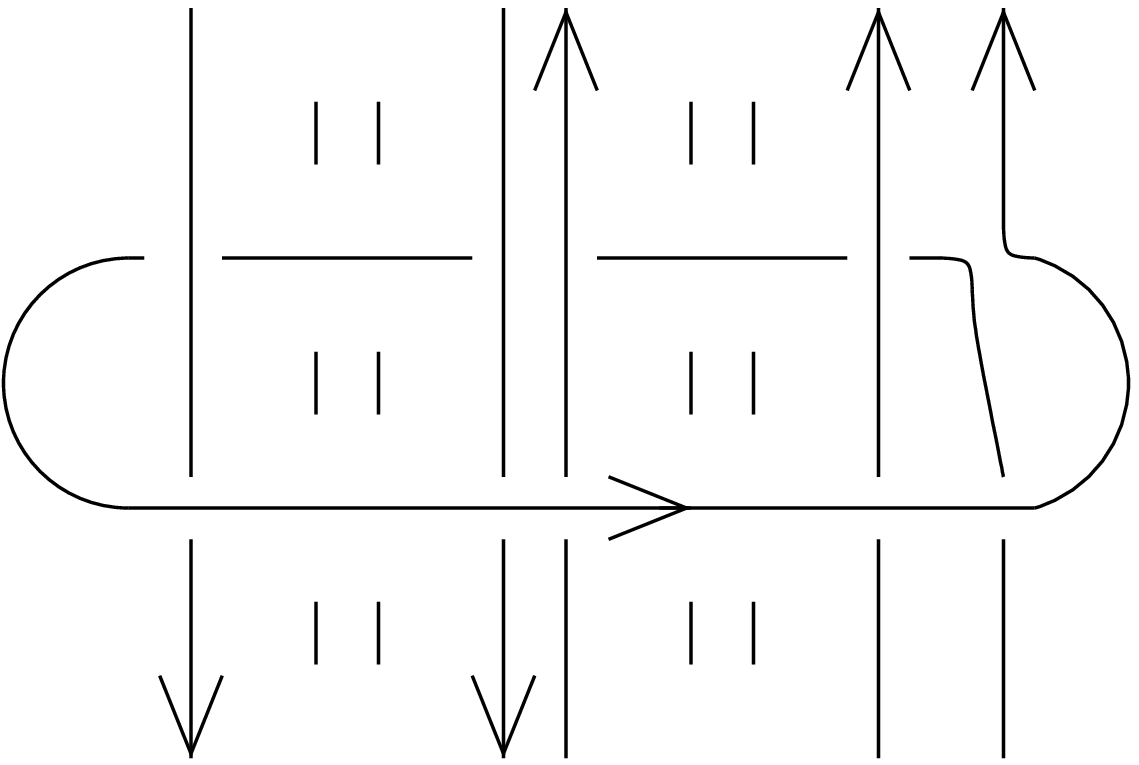}}}\\
&=&T^{(n,p-1)}\otimes 1^{(1)}_{(+)}+zv^{-1}A(p).
\end{eqnarray*} Repeated application of the skein relation in this way will
clearly yield:
\begin{eqnarray} T^{(n,p)}&=&T^{(n,0)}\otimes
1^{(p)}_{(+)}+zv^{-1}\sum_{j=1}^pA(j)\nonumber \\
         &=&T^{(n)}_{(-)}\otimes 1^{(p)}_{(+)}+zv^{-1}\sum_{j=1}^pA(j).
\end{eqnarray}

Now observe, similar to a result in \cite{Morton01}, we can find:
\begin{eqnarray}
{\raisebox{-16pt}{\includegraphics[height=35pt]{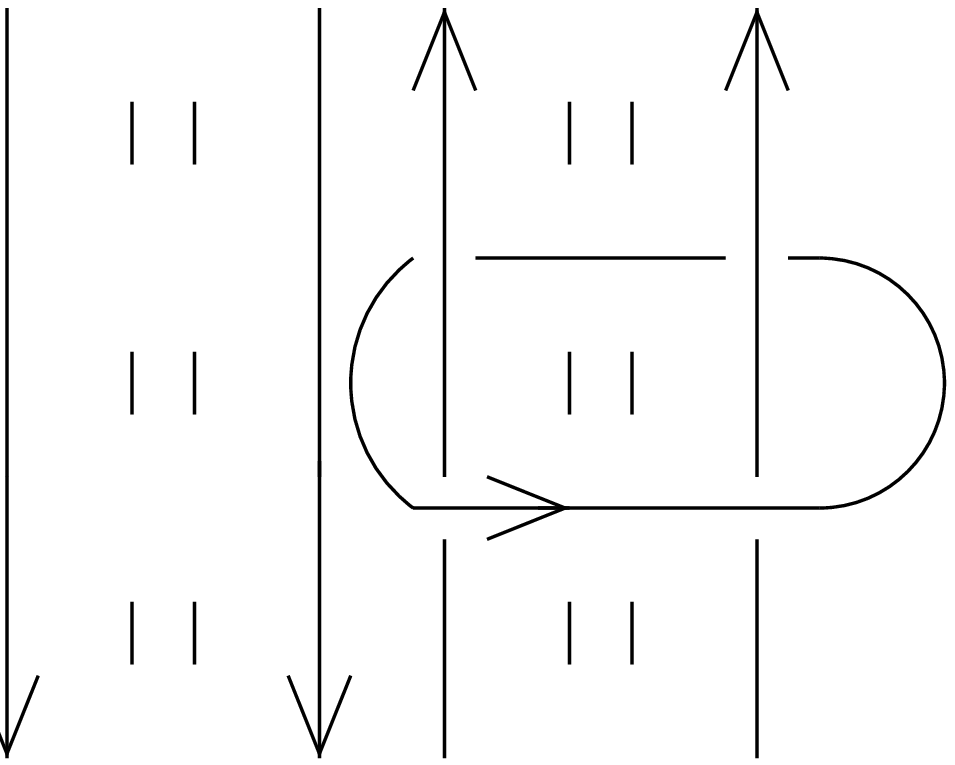}}}&=&\delta\,{\raisebox{-16pt}{\includegraphics[height=35pt]{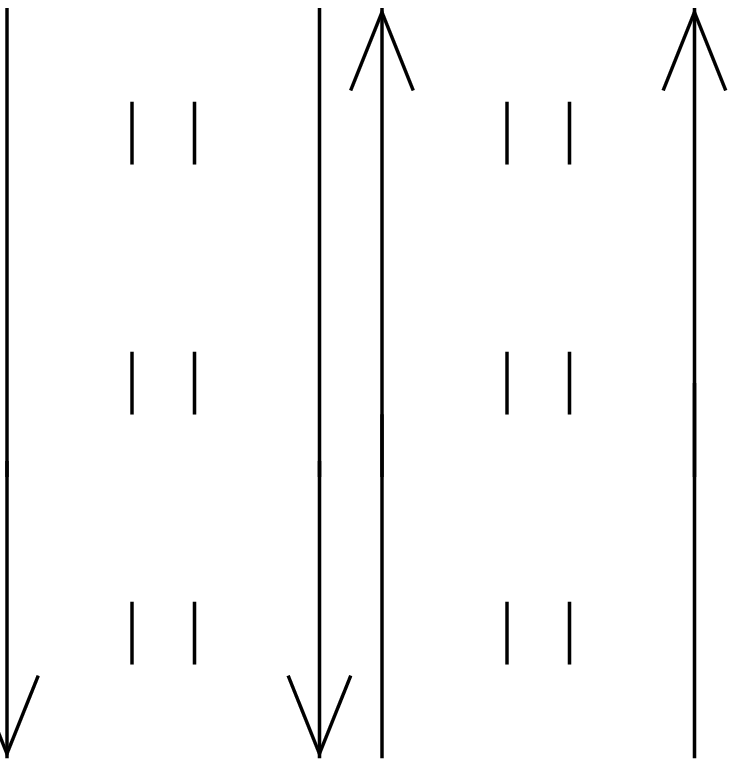}}}+zv^{-1}\sum_{j=1}^pT(j).
\end{eqnarray}

Combining equations $1$ and $2$, we see that we are only left to show that:\[
zv^{-1}\sum_{j=1}^pA(j)=zv^{-1}\sum_{j=1}^pT(j)+w, \] for $w\in M_{n,p}^{(1)}$.

Let $w=\sum_{j=1}^pw(j)$. We must now show that for each $j$, with $1\le j\le
p$, there exists a $w(j)$ such that:
\[ zv^{-1}A(j)=zv^{-1}T(j)+w(j).\] 

Now,
\begin{eqnarray*}
zv^{-1}A(j)&=&zv^{-1}\left\{{\raisebox{-16pt}{\includegraphics[height=35pt]{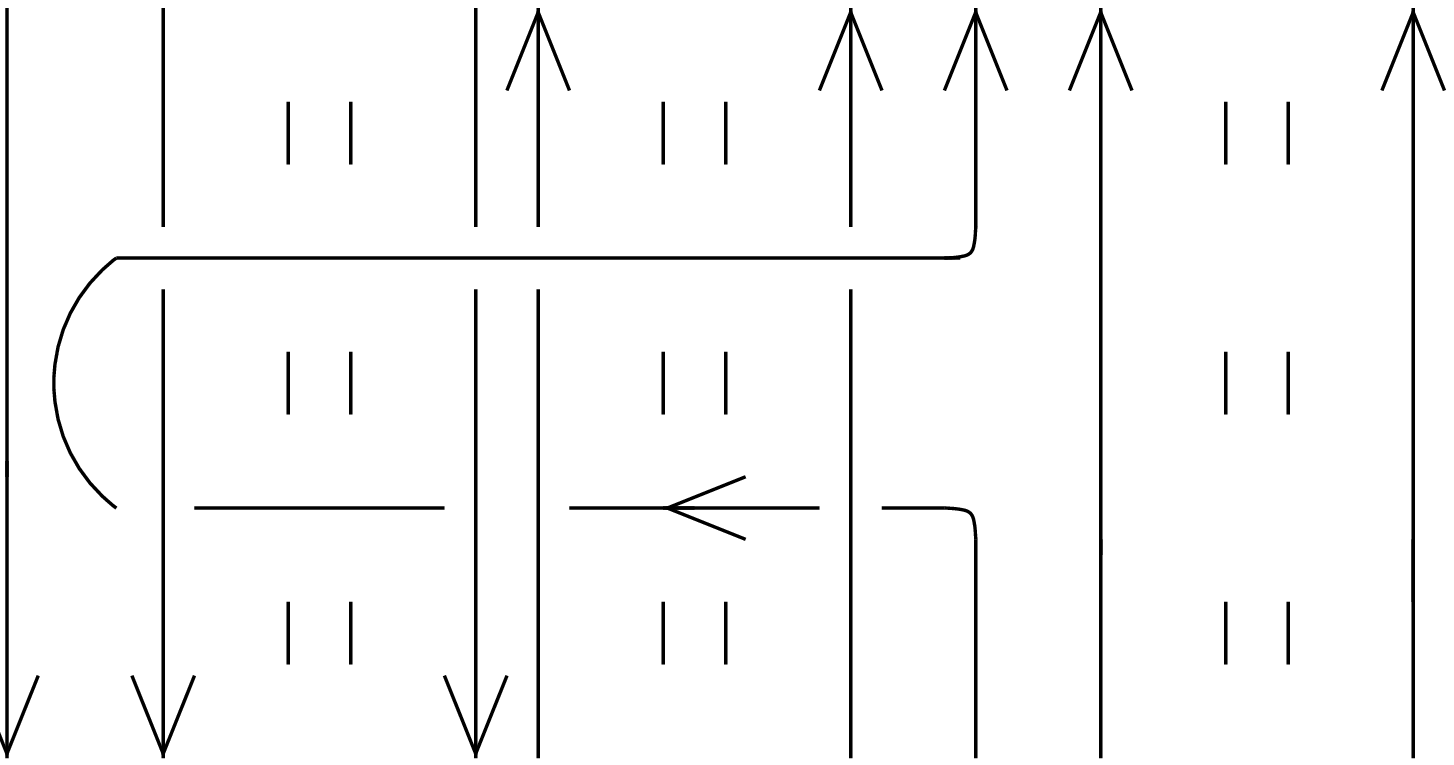}}}-z\,{\raisebox{-16pt}{\includegraphics[height=35pt]{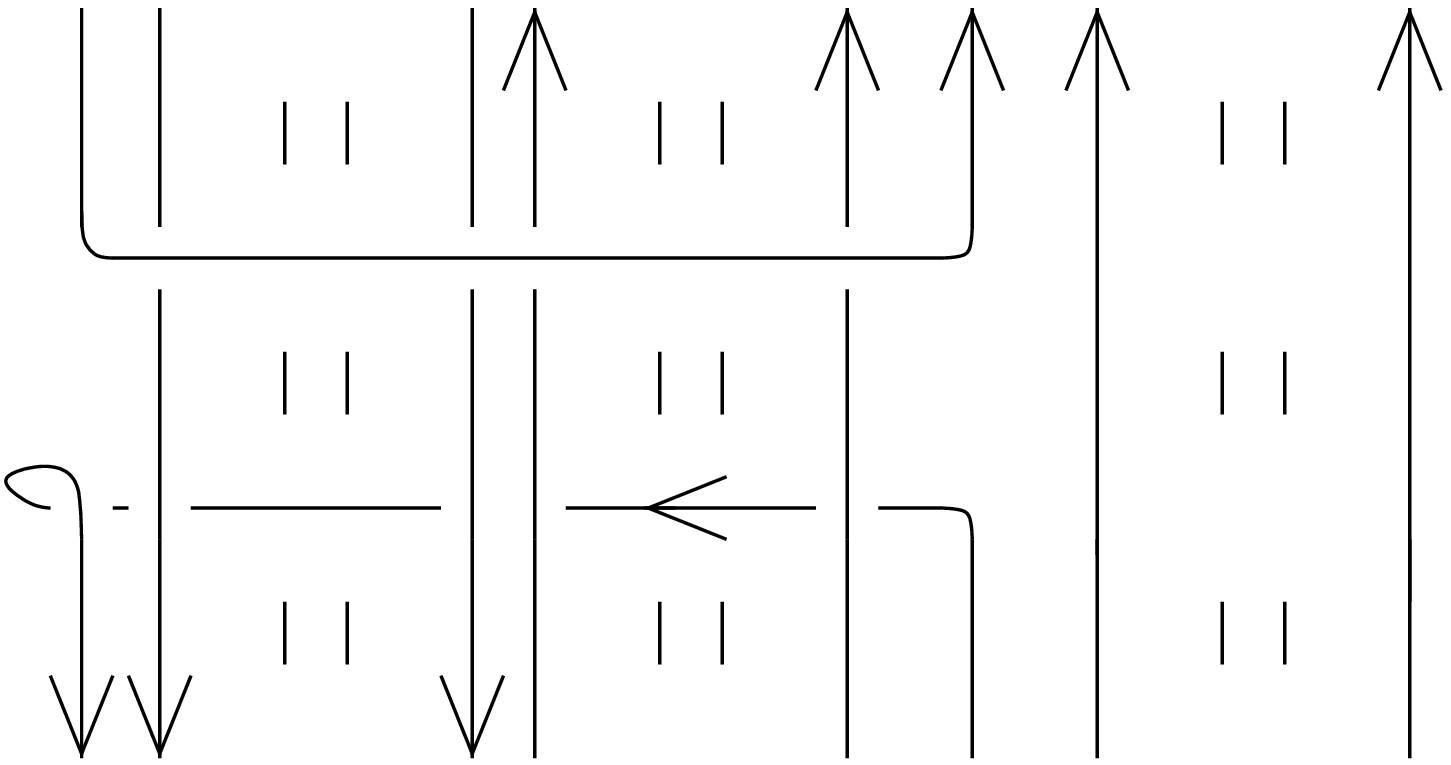}}}\right\}\\
&=&zv^{-1}\left\{{\raisebox{-16pt}{\includegraphics[height=35pt]{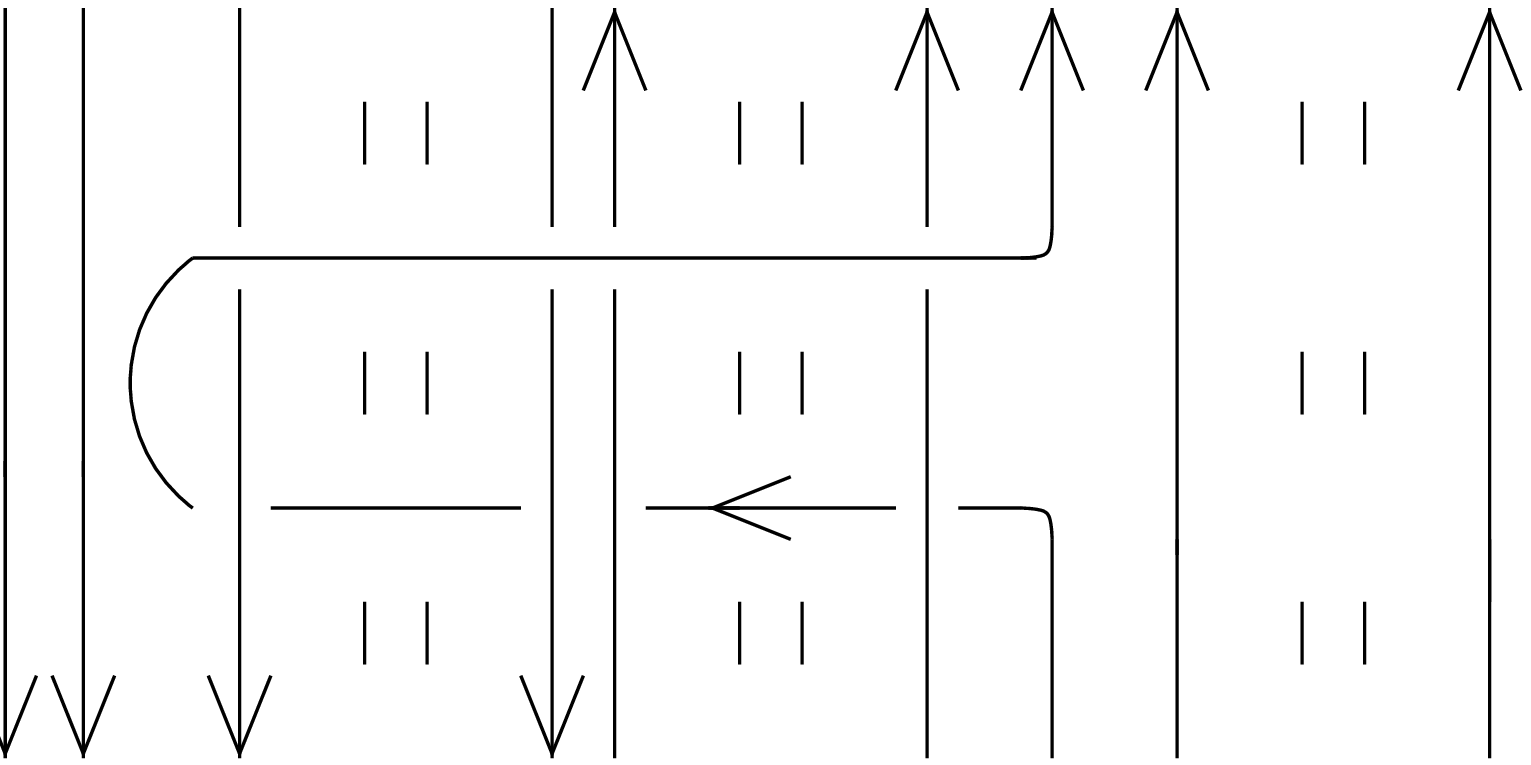}}}-z\,{\raisebox{-16pt}{\includegraphics[height=35pt]{skein5}}}-z\,{\raisebox{-16pt}{\includegraphics[height=35pt]{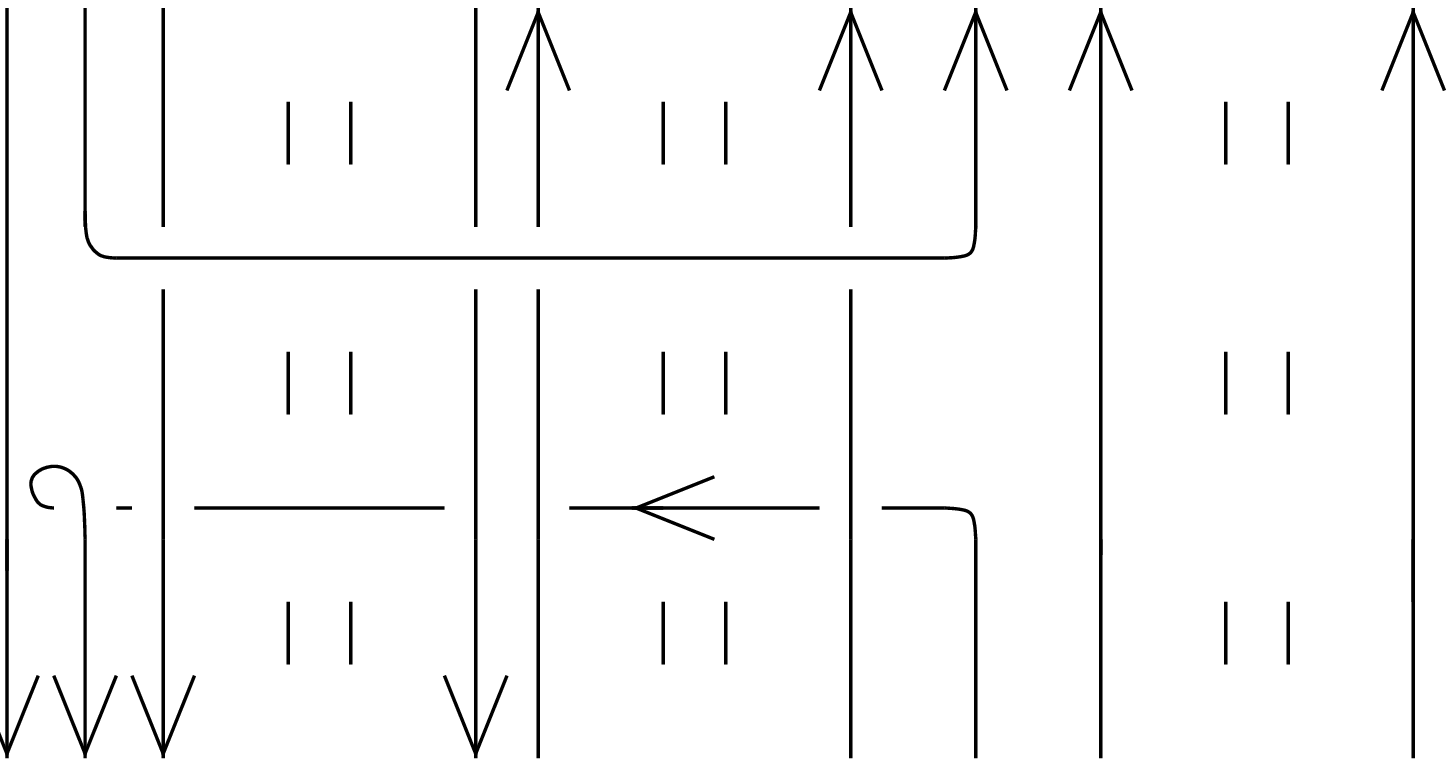}}}\right\}\\
      &=&\cdots\text{ (repeating application of the skein
relation)}\phantom{{\raisebox{-16pt}{\includegraphics[height=35pt]{skein7}}}} \\
&=&zv^{-1}T(j)+z^2v^{-1}\left\{-{\raisebox{-16pt}{\includegraphics[height=35pt]{skein5}}}-\cdots-{\raisebox{-16pt}{\includegraphics[height=35pt]{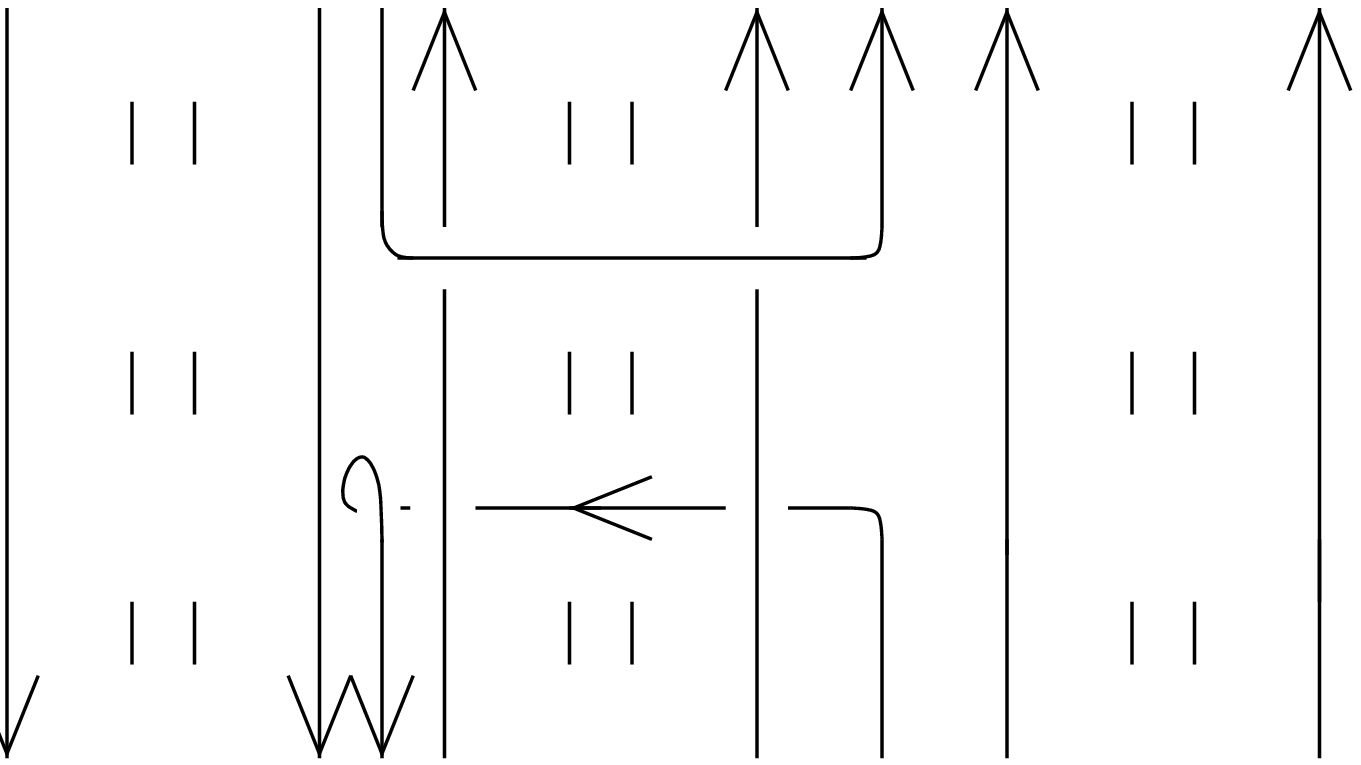}}}\right\}\\
      &=&zv^{-1}T(j)+w(j).
\end{eqnarray*} with $w(j)\in M_{n,p}^{(1)}$.

The result follows.
\end{proof}

We can find an obvious set of quasi-idempotent elements in $M_{n,p}$ given by
$e_{\lambda,\mu}^\prime:=e_\lambda^{(-)}\otimes e_\mu^{(+)}$ formed by the
juxtaposition of the Gyoja-Aiston idempotents with appropriate orientations and
$|\lambda|=n$ and $|\mu|=p$.  There are then $\pi(n)\times\pi(p)$ of these.

We can then use the information in the previous section combined with
Lemma~\ref{thm:decomp} to prove the following proposition.

\begin{prop}\label{thm:tlm}
\begin{eqnarray*}T^{(n,p)}e_{\lambda,\mu}^\prime&=&t_{\lambda,\mu}e_{\lambda,\mu}^\prime+we_{\lambda,\mu}^{\prime}
\\
\text{and}\phantom{--}\bar{T}^{(n,p)}e_{\lambda,\mu}^\prime&=&\bar{t}_{\lambda,\mu}e_{\lambda,\mu}^\prime+w^{\prime}e_{\lambda,\mu}^{\prime},
\end{eqnarray*} where,
\begin{eqnarray*}
t_{\lambda,\mu}&=&(s-s^{-1})\left(-v\sum_{\substack{\text{cells}\\ \text{in
}\lambda}}s^{-2(\text{content})}+v^{-1}\sum_{\substack{\text{cells}\\ \text{in
}\mu}}s^{2(\text{content})}\right)+\delta \\
\text{and}\\
\bar{t}_{\lambda,\mu}&=&(s-s^{-1})\left(v^{-1}\sum_{\substack{\text{cells}\\
\text{in }\lambda}}s^{2(\text{content})}-v\sum_{\substack{\text{cells}\\ \text{in
}\mu}}s^{-2(\text{content})}\right)+\delta.
\end{eqnarray*}
\end{prop} Here we had fixed $|\lambda|$ and $|\mu|$ with values $n$ and $p$
respectively.  In fact, we find that $t_{\lambda,\mu}$ and
$\bar{t}_{\lambda,\mu}$ have the following property:
\begin{lem}\label{thm:distinct} As $\lambda$ and $\mu$ vary over all choices of
Young diagram, the values of $t_{\lambda,\mu}$ are all distinct; as are the
values of $\bar{t}_{\lambda,\mu}$.
\end{lem}
\begin{rem} An equivalent way of stating Lemma~\ref{thm:distinct} is that
\emph{if $t_{\lambda,\mu}=t_{\lambda^\prime,\mu^\prime}$ then
$\lambda=\lambda^\prime$ and $\mu=\mu^\prime$} (similarly for the
$\bar{t}_{\lambda,\mu}$).
\end{rem}
\begin{proof} (of Lemma~\ref{thm:distinct}) We prove the first part of the lemma and note that the second part
follows immediately due to the observation that
$\bar{t}_{\lambda,\mu}=t_{\mu,\lambda}$.

Given $f(s,v)=t_{\lambda,\mu}$ we now show how to recover the Young diagrams
$\lambda$ and $\mu$.

From the formula for $t_{\lambda,\mu}$ in Lemma~\ref{thm:tlm} we see that
$f(s,v)-\delta$ is a Laurent polynomial in $s$ and $v$, and must be of the form:
\[ (s-s^{-1})(-vP(s)+v^{-1}Q(s)).
\]

Now consider $P(s)$ and $Q(s)$ individually.  It is clear that these are also
Laurent polynomials, this time only in the variable $s$.  We have
\begin{eqnarray*} P(s)&=&\sum a_is^{-2i}\\
\text{and}\phantom{--}Q(s)&=&\sum b_js^{2j},
\end{eqnarray*} where $a_i$ is the number of cells in $\lambda$ with content $i$,
and similarly, $b_j$ is the number of cells in $\mu$ with content $j$.  Hence we
can uniquely construct $\lambda$ and $\mu$.
\end{proof}

Now let us return to the maps $\varphi$ and $\bar{\varphi}$, restricting them to
the skein $\annu$.

\begin{thm}\label{thm:eigen} The $t_{\lambda,\mu}$ and $\bar{t}_{\lambda,\mu}$
are eigenvalues of $\varphi|_{\ann^{(n,p)}}$ and $\bar{\varphi}|_{\ann^{(n,p)}}$
respectively.  Moreover, they occur with multiplicity $1$.
\end{thm}
\begin{proof} We prove the result for the $t_{\lambda,\mu}$ with an identical
argument proving the result for the $\bar{t}_{\lambda,\mu}$.

Fix an integer $k$ such that $k=p-n$ and $k\ge 0$ (in other words $p\ge n$ ---
the case for $p<n$ is identical).  Write $\annu$ as $\ann^{(n,k+n)}$ and do
induction on $n$.

For $n=0$ we have that $\ann^{(0,k)}\cong\ann^{(k)}$.  Now for $|\lambda|=0$ and
$|\mu|=k$ we have that $t_{\lambda,\mu}=t_\mu$.  Moreover, in the proof of
Theorem~\ref{thm:mort} we saw that the $t_\mu$ with $|\mu|=k$ are eigenvalues of
$\varphi|_{\ann^{(k)}}$.  Now since
$\ann^{(k)}\cong\ann^{(0,k)}\subset\ann^{(n,k+n)}$ for all $n$, the $t_\mu$ are
also eigenvalues of $\varphi|_{\ann^{(n,k+n)}}$.

Now assume that for $|\lambda|<n$ and $|\mu|<k+n$ the $t_{\lambda,\mu}$ are
eigenvalues of $\varphi|_{\ann^{(|\lambda|,|\mu|)}}$.  Since
$\ann^{(|\lambda|,|\mu|)}\subset\ann^{(n,k+n)}$ the $t_{\lambda,\mu}$ are also
eigenvalues of $\varphi|_{\ann^{(n,k+n)}}$.

Consider the $t_{\lambda,\mu}$ with $|\lambda|=n$ and $|\mu|=k+n$.  By the
inductive hypothesis these $t_{\lambda,\mu}$ are not eigenvalues of
$\varphi|_{\ann^{(n-1,k+n-1)}}$ since we have $\pi(n-1,k+n-1)$ eigenvalues and 
$\ann^{(n-1,k+n-1)}$ is spanned by $\pi(n-1,k+n-1)$ elements and by
Lemma~\ref{thm:distinct} we have that if
$t_{\lambda,\mu}=t_{\lambda^\prime,\mu^\prime}$ then $\lambda=\lambda^\prime$ and
$\mu=\mu^\prime$.

Define elements $Q^\prime_{\lambda,\mu}:=Q^{(-)}_\lambda\cdot
Q^{(+)}_\mu(=\wedge(e^\prime_{\lambda,\mu}))$ with $|\lambda|=n$ and
$|\mu|=k+n$.  Clearly $Q^\prime_{\lambda,\mu}\in\ann^{(n,k+n)}$.

Now by Lemma~\ref{thm:tlm},
\[
\varphi|_{\ann^{(n,k+n)}}(Q^\prime_{\lambda,\mu})=t_{\lambda,\mu}Q^\prime_{\lambda,\mu}+w^\prime
\] where $w^\prime\in\ann^{(n-1,k+n-1)}$.

We can find a $v\in\ann^{(n-1,k+n-1)}$ such that
$(\varphi|_{\ann^{(n,k+n)}}-t_{\lambda,\mu}I)(v)=w^\prime$.

Now consider $Q^\prime_{\lambda,\mu}-v$.  This is clearly non-zero.  We find:
\begin{eqnarray*}
\varphi|_{\ann^{(n,k+n)}}(Q^\prime_{\lambda,\mu}-v)&=&\varphi|_{\ann^{(n,k+n)}}(Q^\prime_{\lambda,\mu})-\varphi|_{\ann^{(n,k+n)}}(v)+t_{\lambda,\mu}v-t_{\lambda,\mu}v
\\ &=&\varphi|_{\ann^{(n,k+n)}}(Q^\prime_{\lambda,\mu})-w^\prime-t_{\lambda,\mu}v
\\ &=&t_{\lambda,\mu}Q^\prime_{\lambda,\mu}+w^\prime-w^\prime-t_{\lambda,\mu}v \\
&=&t_{\lambda,\mu}(Q^\prime_{\lambda,\mu}-v).
\end{eqnarray*} Hence such $t_{\lambda,\mu}$ are eigenvalues of
$\varphi|_{\ann^{(n,k+n)}}$.

Hence by induction, we have that the $t_{\lambda,\mu}$, with $|\lambda|\le n$,
$|\mu|\le p$ and $|\lambda|-|\mu|=n-p$, are eigenvalues of
$\varphi|_{\ann^{(n,p)}}$.

Moreover, we have found at least $\pi(n,p)$ eigenvalues for
$\varphi|_{\ann^{(n,p)}}$.  But $\ann^{(n,p)}$ is known to be spanned by
$\pi(n,p)$ elements, so $\varphi|_{\ann^{(n,p)}}$ has at most $\pi(n,p)$ different
eigenvalues.  Hence it has exactly $\pi(n,p)$ eigenvalues each with multiplicity
one.
\end{proof}

We now state two useful corollaries.

\begin{cor} There is a basis of $\annu$ given by:
\[
\{Q_{\lambda,\mu}:|\lambda|\le n,|\mu|\le p,|\lambda|-|\mu|=n-p\}
\] such that:
\[
\varphi(Q_{\lambda,\mu})=t_{\lambda,\mu}Q_{\lambda,\mu}\phantom{--}\text{and}\phantom{--}\bar{\varphi}(Q_{\lambda,\mu})=\bar{t}_{\lambda,\mu}Q_{\lambda,\mu}.
\]
\end{cor}
\begin{cor} Every eigenvector of $\varphi$ and $\bar{\varphi}$ is a multiple of
one such basis element.
\end{cor}
\begin{rem} The eigenvalues $t_{\lambda,\mu}$ and $\bar{t}_{\lambda,\mu}$
correspond to the eigenvalues of the matrix $M$ in equation (1.1) of
\cite{Chan00}, found there only for $1\le k_1+k_2\le 5$ and $k_2\le k_1$.  Chan
uses the Homfly polynomial based on parameters $l$ and $m$, which are variants of
$v$ and $z$.  The numbers $\sqrt{m^2-4}$ in Chan's eigenvalues $\rho_i$ and
$\rho^*_i$ correspond to the parameter $s$ here with $z=s-s^{-1}$, which features
strongly in our eigenvalues $t_{\lambda,\mu}$ and $\bar{t}_{\lambda,\mu}$.  Our
use of $s$ is the feature which allows us to give simple formulae for the
Gyoja-Aiston elements $Q_\lambda$ and to extend in principle to $Q_{\lambda,\mu}$.
\end{rem}

Unlike the Gyoja-Aiston elements $Q_\lambda$ which are known and have been
well-studied, their generalisations the $Q_{\lambda,\mu}$ described in the above
Corollary are not well-understood.  We shall show in the following section how
they can be found explicitly.

\section{The Homfly Polynomials of Some Generalized\break Hopf Links} Here we apply the
techniques described above to show how computation of the Homfly polynomial of
some generalized Hopf links is possible.

\subsection{The Homfly polynomial of $H(k_1,k_2;n,0)\,(=H^*(k_1,k_2;0,n))$}
Consider $H(k_1,k_2;n,0)$ in the skein of the annulus.  Then we have
\[ H(k_1,k_2;n,0)=\varphi^{k_1}(\bar{\varphi}^{k_2}(A_1^n).
\]

Now since the maps $\varphi$ and $\bar{\varphi}$ are linear maps, we know that
for the $Q_\lambda$,
\[
\varphi^{k_1}(\bar{\varphi}^{k_2}(Q_\lambda))=t_\lambda^{k_1}\bar{t}_\lambda^{\,k_2}Q_\lambda.
\]

Also, since the $Q_\lambda$ are a basis or the skein $\ann^{(n)}$, we have
\[ A_1^n=\sum_{|\lambda|=n}d_\lambda Q_\lambda
\] for constants $d_\lambda$.  The $d_\lambda$ can be calculated by several
means, for example by counting the number of standard tableaux of shape $\lambda$.

Therefore,
\begin{eqnarray*}
H(k_1,k_2;n,0)&=&\sum_{|\lambda|=n}d_\lambda\varphi^{k_1}(\bar{\varphi}^{k_2}(Q_\lambda))
\\ &=&\sum_{|\lambda|=n}d_\lambda t_\lambda^{k_1}\bar{t}_\lambda^{k_2}Q_\lambda.
\end{eqnarray*} So evaluating in the plane (using the work of \cite{AM98}), we
find
\[ P(H(k_1,k_2;n,0))=\sum_{|\lambda|=n}d_\lambda
t_\lambda^{k_1}\bar{t}_\lambda^{k_2}\left(\prod_{(i,j)\in\lambda}\frac{v^{-1}s^{j-i}-vs^{i-j}}{s^{\text{hl}(i,j)}-s^{-\text{hl}(i,j)}}\right),
\] where $\text{hl}(i,j)$ is the hook-length of the cell $(i,j)$, in row $i$ and
column $j$.

\subsection{The Homfly polynomial of $H(k_1,k_2;n_1,n_2)$} Consider, in a similar
way to above, $H(k_1,k_2;n_1,n_2)$ as an element of the skein $\ann$.  Then we
have
\[ H(k_1,k_2;n_1,n_2)=\varphi^{k_1}(\bar{\varphi}^{k_2}(A_1^{n_1}A_{-1}^{n_2})).
\] Similar to the restricted case above, we have
\[
\varphi^{k_1}(\bar{\varphi}^{k_2}(Q_{\lambda,\mu}))=t_{\lambda,\mu}^{k_1}\bar{t}_{\lambda,\mu}^{\,k_2}Q_{\lambda,\mu}
\] and
\[ A_1^{n_1}A_{-1}^{n_2}=\sum_{\substack{|\lambda|\le n_2 \\ |\mu|\le n_1 \\
|\lambda|-|\mu|=n_2-n_1}}d_{\lambda,\mu}Q_{\lambda,\mu}
\] for constants $d_{\lambda,\mu}$.  These constants can be calculated in terms
of appropriate $d_\lambda$ and $d_\mu$ (see previous section).
\begin{thm}\label{thm:stem}{\rm\cite{Stem87}}\qua The numbers $d_{\lambda,\mu}$ can be
found from the following formula:
\[ d_{\lambda,\mu}=m!{n_2 \choose m}{n_1 \choose m}d_\lambda d_\mu,
\] where $|\lambda|\le n_2$, $|\mu|\le n_1$ and $m=n_2-|\lambda|=n_1-|\mu|$.
\end{thm} Therefore,
\begin{eqnarray*} H(k_1,k_2;n_1,n_2)&=&\sum_{\substack{|\lambda|\le n_2 \\
|\mu|\le n_1 \\
|\lambda|-|\mu|=n_2-n_1}}d_{\lambda,\mu}\varphi^{k_1}(\bar{\varphi}^{k_2}(Q_{\lambda,\mu}))
\\ &=&\sum_{\substack{|\lambda|\le n_2 \\ |\mu|\le n_1 \\
|\lambda|-|\mu|=n_2-n_1}}d_{\lambda,\mu}t_{\lambda,\mu}^{k_1}\bar{t}_{\lambda,\mu}^{\,k_2}Q_{\lambda,\mu}.
\end{eqnarray*} At present, we do not have a general closed formula for
$P(H(k_1,k_2;n_1,n_2))$ due to lack of information about the elements
$Q_{\lambda,\mu}$.

We can, however, make explicit calculations in individual cases as illustrated by
the following example.
\begin{exmp} Consider $H(k_1,k_2;1,2)\in\ann^{(2,1)}$, as shown in
Figure~\ref{fig:example}.
\begin{figure}[ht!]
\psfrag{3}{$k_1$}
\psfrag{4}{$k_2$}
\begin{center}
\includegraphics[width=0.45\textwidth,height=6cm]{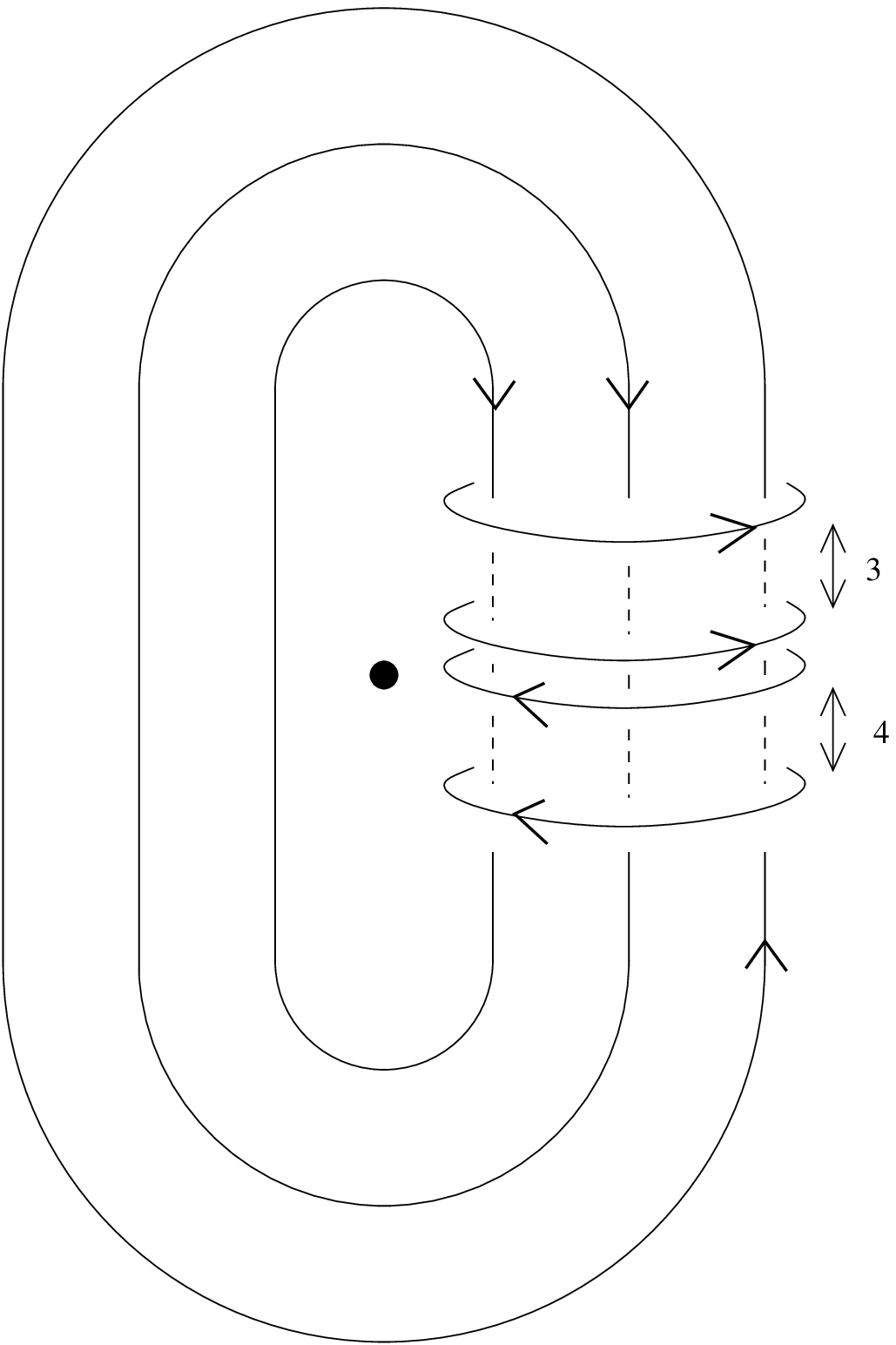}
\caption{The link $H(k_1,k_2;1,2)$ in $\ann$}
\label{fig:example}
\end{center}
\end{figure}

Then
\[ H(k_1,k_2;1,2)=\varphi^{k_1}(\bar{\varphi}^{k_2}(A_1A_{-1}^2)),
\] where, by Theorem~\ref{thm:stem},
\begin{eqnarray}
A_1A_{-1}^2=Q_{\raisebox{-3pt}{{\includegraphics[height=5pt]{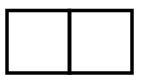}},{\includegraphics[height=5pt]{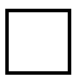}}}}+2Q_{\raisebox{-3pt}{{\includegraphics[height=5pt]{1}},$\scriptstyle{\emptyset}$}}+Q_{\raisebox{-3pt}{\raisebox{-3pt}{{\includegraphics[height=10pt]{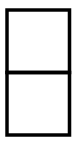}}},{\includegraphics[height=5pt]{1}}}}.
\end{eqnarray} However, we can also find, by using powers of trivial Gyoja-Aiston
elements $Q_{\raisebox{-3pt}{{\includegraphics[height=5pt]{1}}}}$, with
appropriate orientation, that
\[ A_1A_{-1}^2=(Q_{\raisebox{-3pt}{{\includegraphics[height=5pt]{1}}}}^{(-)})^2
Q_{\raisebox{-3pt}{{\includegraphics[height=5pt]{1}}}}^{(+)}.
\] Moreover, these elements are known to satisfy the Littlewood-Richardson rule
for multiplication of Young diagrams (\cite{AistonPhD}), so
\begin{eqnarray}
A_1A_{-1}^2&=&(Q_{\raisebox{-3pt}{{\includegraphics[height=5pt]{2}}}}^{(-)}+Q_{\raisebox{-3pt}{\raisebox{-3pt}{{\includegraphics[height=10pt]{11}}}}}^{(-)})Q_{\raisebox{-3pt}{{\includegraphics[height=5pt]{1}}}}^{(+)}
\nonumber \\
&=&Q_{\raisebox{-3pt}{{\includegraphics[height=5pt]{2}}}}^{(-)}Q_{\raisebox{-3pt}{{\includegraphics[height=5pt]{1}}}}^{(+)}+Q_{\raisebox{-3pt}{{\includegraphics[height=10pt]{11}}}}^{(-)}Q_{\raisebox{-3pt}{{\includegraphics[height=5pt]{1}}}}^{(+)}
\nonumber \\
&=&Q_{\raisebox{-3pt}{{\includegraphics[height=5pt]{2}},{\includegraphics[height=5pt]{1}}}}^{\prime}+Q_{\raisebox{-3pt}{\raisebox{-3pt}{{\includegraphics[height=10pt]{11}}},{\includegraphics[height=5pt]{1}}}}^{\prime}.
\end{eqnarray} Now combining equations (3) and (4) with the observation that 
\[Q_{\raisebox{-3pt}{{\includegraphics[height=5pt]{1}},$\scriptstyle{\emptyset}$}}=Q_{\raisebox{-3pt}{{\includegraphics[height=5pt]{1}},$\scriptstyle{\emptyset}$}}^{\prime}=Q_{\raisebox{-3pt}{{\includegraphics[height=5pt]{1}}}}^{(-)}Q_{\raisebox{-3pt}{$\scriptstyle{\emptyset}$}}^{(+)}\]
and assuming symmetry under conjugation of Young diagrams, we have:
\begin{eqnarray*}
Q_{\raisebox{-3pt}{{\includegraphics[height=5pt]{2}},{\includegraphics[height=5pt]{1}}}}&=&Q_{\raisebox{-3pt}{{\includegraphics[height=5pt]{2}},{\includegraphics[height=5pt]{1}}}}^{\prime}-Q_{\raisebox{-3pt}{{\includegraphics[height=5pt]{1}},$\scriptstyle{\emptyset}$}}^{\prime},\\
\text{and
}Q_{\raisebox{-3pt}{\raisebox{-3pt}{{\includegraphics[height=10pt]{11}}},{\includegraphics[height=5pt]{1}}}}&=&Q_{\raisebox{-3pt}{\raisebox{-3pt}{{\includegraphics[height=10pt]{11}}},{\includegraphics[height=5pt]{1}}}}^{\prime}-Q_{\raisebox{-3pt}{{\includegraphics[height=5pt]{1}},$\scriptstyle{\emptyset}$}}^{\prime}.
\end{eqnarray*} Hence, evaluating in the plane, we find,
\begin{eqnarray}
P(H(k_1,k_2;1,2))&=&P(\varphi^{k_1}(\bar{\varphi}^{k_2}(A_1A_{-1}^2))) \nonumber
\\
&=&t^{k_1}_{\raisebox{-3pt}{{\includegraphics[height=5pt]{2}},{\includegraphics[height=5pt]{1}}}}\bar{t}^{k_2}_{\raisebox{-3pt}{{\includegraphics[height=5pt]{2}},{\includegraphics[height=5pt]{1}}}}P(Q_{\raisebox{-3pt}{{\includegraphics[height=5pt]{2}},{\includegraphics[height=5pt]{1}}}})
\nonumber \\
& &+2t^{k_1}_{\raisebox{-3pt}{{\includegraphics[height=5pt]{1}},$\scriptstyle{\emptyset}$}}\bar{t}^{k_2}_{\raisebox{-3pt}{{\includegraphics[height=5pt]{1}},$\scriptstyle{\emptyset}$}}P(Q_{\raisebox{-3pt}{{\includegraphics[height=5pt]{1}},$\scriptstyle{\emptyset}$}})+t^{k_1}_{\raisebox{-3pt}{\raisebox{-3pt}{{\includegraphics[height=10pt]{11}}},{\includegraphics[height=5pt]{1}}}}\bar{t}^{k_2}_{\raisebox{-3pt}{\raisebox{-3pt}{{\includegraphics[height=10pt]{11}}},{\includegraphics[height=5pt]{1}}}}P(Q_{\raisebox{-3pt}{\raisebox{-3pt}{{\includegraphics[height=10pt]{11}}},{\includegraphics[height=5pt]{1}}}})
\nonumber \\
&=&t^{k_1}_{\raisebox{-3pt}{{\includegraphics[height=5pt]{2}},{\includegraphics[height=5pt]{1}}}}\bar{t}^{k_2}_{\raisebox{-3pt}{{\includegraphics[height=5pt]{2}},{\includegraphics[height=5pt]{1}}}}(P(Q_{\raisebox{-3pt}{{\includegraphics[height=5pt]{2}},{\includegraphics[height=5pt]{1}}}}^{\prime})-P(Q_{\raisebox{-3pt}{{\includegraphics[height=5pt]{1}},$\scriptstyle{\emptyset}$}}^{\prime}))
\nonumber \\
& &+2t^{k_1}_{\raisebox{-3pt}{{\includegraphics[height=5pt]{1}},$\scriptstyle{\emptyset}$}}\bar{t}^{k_2}_{\raisebox{-3pt}{{\includegraphics[height=5pt]{1}},$\scriptstyle{\emptyset}$}}P(Q_{\raisebox{-3pt}{{\includegraphics[height=5pt]{1}},$\scriptstyle{\emptyset}$}}^{\prime})+t^{k_1}_{\raisebox{-3pt}{\raisebox{-3pt}{{\includegraphics[height=10pt]{11}}},{\includegraphics[height=5pt]{1}}}}\bar{t}^{k_2}_{\raisebox{-3pt}{\raisebox{-3pt}{{\includegraphics[height=10pt]{11}}},{\includegraphics[height=5pt]{1}}}}(P(Q^{\prime}_{\raisebox{-3pt}{\raisebox{-3pt}{{\includegraphics[height=10pt]{11}}},{\includegraphics[height=5pt]{1}}}}-P(Q_{\raisebox{-3pt}{{\includegraphics[height=5pt]{1}},$\scriptstyle{\emptyset}$}}^{\prime}))
\nonumber \\
&=&t^{k_1}_{\raisebox{-3pt}{{\includegraphics[height=5pt]{2}},{\includegraphics[height=5pt]{1}}}}\bar{t}^{k_2}_{\raisebox{-3pt}{{\includegraphics[height=5pt]{2}},{\includegraphics[height=5pt]{1}}}}P(Q_{\raisebox{-3pt}{{\includegraphics[height=5pt]{2}},{\includegraphics[height=5pt]{1}}}}^{\prime})
\nonumber \\
                 &
&+(2t^{k_1}_{\raisebox{-3pt}{{\includegraphics[height=5pt]{1}},$\scriptstyle{\emptyset}$}}\bar{t}^{k_2}_{\raisebox{-3pt}{{\includegraphics[height=5pt]{1}},$\scriptstyle{\emptyset}$}}-t^{k_1}_{\raisebox{-3pt}{{\includegraphics[height=5pt]{2}},{\includegraphics[height=5pt]{1}}}}\bar{t}^{k_2}_{\raisebox{-3pt}{{\includegraphics[height=5pt]{2}},{\includegraphics[height=5pt]{1}}}}-t^{k_1}_{\raisebox{-3pt}{\raisebox{-3pt}{{\includegraphics[height=10pt]{11}}},{\includegraphics[height=5pt]{1}}}}\bar{t}^{k_2}_{\raisebox{-3pt}{\raisebox{-3pt}{{\includegraphics[height=10pt]{11}}},{\includegraphics[height=5pt]{1}}}})P(Q_{\raisebox{-3pt}{{\includegraphics[height=5pt]{1}},$\scriptstyle{\emptyset}$}}^{\prime})
\\
                 &
&+t^{k_1}_{\raisebox{-3pt}{\raisebox{-3pt}{{\includegraphics[height=10pt]{11}}},{\includegraphics[height=5pt]{1}}}}\bar{t}^{k_2}_{\raisebox{-3pt}{\raisebox{-3pt}{{\includegraphics[height=10pt]{11}}},{\includegraphics[height=5pt]{1}}}}P(Q_{\raisebox{-3pt}{\raisebox{-3pt}{{\includegraphics[height=10pt]{11}}},{\includegraphics[height=5pt]{1}}}}^{\prime})
\nonumber
\end{eqnarray} From the definition of the $Q_{\lambda,\mu}^{\prime}$, we can now
use the results in \cite{AM98} to find
$P(Q_{\raisebox{-3pt}{{\includegraphics[height=5pt]{1}},$\scriptstyle{\emptyset}$}}^{\prime})$,
$P(Q_{\raisebox{-3pt}{{\includegraphics[height=5pt]{2}},{\includegraphics[height=5pt]{1}}}}^{\prime})$
and
$P(Q_{\raisebox{-3pt}{\raisebox{-3pt}{{\includegraphics[height=10pt]{11}}},{\includegraphics[height=5pt]{1}}}}^{\prime})$. 
We have:
\begin{eqnarray*}
P(Q_{\raisebox{-3pt}{{\includegraphics[height=5pt]{1}},$\scriptstyle{\emptyset}$}}^{\prime})&=&\frac{v^{-1}-v}{s-s^{-1}},
\\
P(Q_{\raisebox{-3pt}{{\includegraphics[height=5pt]{2}},{\includegraphics[height=5pt]{1}}}}^{\prime})&=&\left(\frac{v^{-1}-v}{s^2-s^{-2}}\right)\left(\frac{v^{-1}s-vs^{-1}}{s-s^{-1}}\right)\left(\frac{v^{-1}-v}{s-s^{-1}}\right),
\\
\text{and
}P(Q_{\raisebox{-3pt}{\raisebox{-3pt}{{\includegraphics[height=10pt]{11}}},{\includegraphics[height=5pt]{1}}}}^{\prime})&=&\left(\frac{v^{-1}-v}{s^2-s^{-2}}\right)\left(\frac{v^{-1}s^{-1}-vs}{s-s^{-1}}\right)\left(\frac{v^{-1}-v}{s-s^{-1}}\right).
\end{eqnarray*} Then using Proposition~\ref{thm:tlm} we find:
\begin{eqnarray*}
t_{\raisebox{-3pt}{{\includegraphics[height=5pt]{1}},$\scriptstyle{\emptyset}$}}&=&-v(s-s^{-1})+\delta,
\\
t_{\raisebox{-3pt}{{\includegraphics[height=5pt]{2}},{\includegraphics[height=5pt]{1}}}}&=&(s-s^{-1})(-v(1+s^{-2})+v^{-1})+\delta,
\\
t_{\raisebox{-3pt}{\raisebox{-3pt}{{\includegraphics[height=10pt]{11}}},{\includegraphics[height=5pt]{1}}}}&=&(s-s^{-1})(-v(1+s^2)+v^{-1})+\delta,
\end{eqnarray*} and
\begin{eqnarray*}
\bar{t}_{\raisebox{-3pt}{{\includegraphics[height=5pt]{1}},$\scriptstyle{\emptyset}$}}&=&v^{-1}(s-s^{-1})+\delta,
\\
\bar{t}_{\raisebox{-3pt}{{\includegraphics[height=5pt]{2}},{\includegraphics[height=5pt]{1}}}}&=&(s-s^{-1})(v^{-1}(1+s^{2})-v)+\delta,
\\
\bar{t}_{\raisebox{-3pt}{\raisebox{-3pt}{{\includegraphics[height=10pt]{11}}},{\includegraphics[height=5pt]{1}}}}&=&(s-s^{-1})(v^{-1}(1+s^{-2})-v)+\delta.
\end{eqnarray*}

Substitution of these values into equation (5) then gives $P(H(k_1,k_2;1,2))$
immediately.
\end{exmp}

\subsection{A final remark} We can in principle write any given element of the
skein $X\in\ann$ as a linear combination of the basis elements
$Q_{\lambda,\mu}$.  Therefore, one can find $\varphi(X)$ and $\bar{\varphi}(X)$,
and hence readily evaluate the Homfly polynomial of
$H(k_1,k_2;X):=H_+(X,A_1^{k_1}A_{-1}^{k_2})$.  The special case
$X=A_1^{n_1}A_{-1}^{n_2}$ gives $H(k_1,k_2;n_1,n_2)$.

% Thanks
\rk{Acknowledgments}

The second author was supported by EPSRC grant 99801479.

\Addresses\recd

\end{document}